%&amstex
\input amstex
\input amsppt.sty

\define\su{\sum\limits}
\define\pro{\prod\limits}
\define\inte{\int\limits}
\define\siz{\sigma_n^{(zz')}}
\define\rz{\rho_n^{(zz')}}
\NoRunningHeads
\magnification 1200
\topmatter 
\title 
Point Processes and the\\ Infinite Symmetric Group \\Part II: Higher Correlation Functions
\endtitle
\author
Alexei Borodin
\endauthor
\thanks 
Supported by the Russian Program for Support of Scientific Schools (grant
96-15-96060).
\endthanks
\abstract
We continue the study of the correlation functions for the point
stochastic processes introduced in Part I (G.~Olshanski).  We find an
integral representation of all the correlation 
functions and their explicit expression in terms of
multivariate hypergeometric functions. Then we define a modification
(``lifting'') of the processes which results in a substantial
simplification of the structure of the correlation functions. It
turns out that the ``lifted'' correlation functions are given by a
determinantal formula involving a kernel. The latter has the form
$(A(x)B(y)-B(x)A(y))/(x-y)$, where $A$ and $B$ are certain Whittaker
functions. Such a form for correlation functions is well known in
the random matrix theory and mathematical physics. Finally, we get some
asymptotic formulas for the correlation functions which are employed
in Part III (A.~Borodin and G.~Olshanski).
\endabstract
\address
{\rm Department of Mathematics, The University of Pennsylvania, Philadelphia,
 PA 19104-6395, U.S.A.} {\it E-mail address\/}: {\tt borodine\@math.upenn.edu}
\endaddress
\toc
\head  {} Introduction  \endhead
\head 1. General structure of  the correlation functions  \endhead
\subhead 1.1.  The moment problem \endsubhead
\subhead 1.2. Combinatorics of the Murnaghan-Nakayama rule \endsubhead
\subhead 1.3. A simplification of the moment problem \endsubhead
\subhead 1.4. Contraction, pseudoconvolution, and diagonalization \endsubhead
\subhead 1.5.  The distributions $\siz(T)$ \endsubhead
\subhead 1.6. Construction of $I_d$\endsubhead
\subhead 1.7. Combinatorics of the substantial structures \endsubhead
\head 2. The correlation functions in positive (negative) hyperoctants \endhead
\subhead 2.1. Why positive hyperoctants are simpler \endsubhead
\subhead 2.2. An integral representation \endsubhead 
\subhead 2.3.* Lauricella hypergeometric functions of type B \endsubhead 
\subhead 2.4.* Expression of the correlation functions via $F_B^{[m]}(a,b;c|y)$\endsubhead 
\subhead 2.5. Simplicity of the processes \endsubhead
\head 3. Lifted processes \endhead
\subhead 3.1. Lifting. Lifted correlation functions \endsubhead
\subhead 3.2. General structure of  $\tilde{\rho}_n^{(zz')}$\endsubhead
\subhead 3.3. Lifted correlation functions in positive (negative) hyperoctants \endsubhead
\head 4. Asymptotics at the origin\endhead
\subhead 4.1. Lifted processes \endsubhead
\subhead 4.2.* Analytic continuation of $F_B^{[m]}(a,b;c|y)$\endsubhead
\subhead 4.3.* Asymptotics of the non-lifted correlation functions \endsubhead 
\head  {} References \endhead
\endtoc
\endtopmatter

\document
\head Introduction \endhead

 In this paper we compute the higher correlation functions of a remarkable family of stochastic point processes introduced by G.I Olshanski in [O]. The main result is an explicit integral representation for the correlation functions. We also introduce a modification of the processes  that provides a connection with certain class of known processes.
The present paper is the continuation of [O]. In the next paper [BO] we give some applications of our results.

In general, our point processes have a representation-theoretic origin, see [KOV], [O]. Here, to be concrete, we give a shorter definition of the processes via finite-point approximations, see [O].

Let us consider the {\it Thoma simplex}  $\Omega$ (see [O])
$$
\Omega= \left\{ {\alpha}_{1}\geq{\alpha}_{2}\geq\ldots\geq 0,\quad
{\beta}_{1}\geq{\beta}_{2}\geq\ldots\geq 0;\quad\sum{(\alpha_i+\beta_i)}\leq 1
\right\}
$$
and its face 
$$
\Omega_0= \left\{ {\alpha}_{1}\geq{\alpha}_{2}\geq\ldots\geq 0,\quad
{\beta}_{1}\geq{\beta}_{2}\geq\ldots\geq 0;\quad\sum{(\alpha_i+\beta_i)}= 1
\right\}.
$$
We define a sequence $P^{(n)}_{zz'};\ n=1,2,\ldots$, of  discrete probability measures on $\Omega_0$ depending on two complex parameters $z,z'$ as follows. The support of the $n$th measure consists of all points of the form
 $$
\left(\frac{p_1+1/2}{n}>  \ldots >\frac{p_d+1/2}{n},\quad               
\frac{q_1+1/2}{n}>\ldots >\frac{q_d+1/2}{n}\right)\in \Omega_0,
$$ 
$$
p_1,q_1,\ldots p_d,q_d\in \{0,1,2,\ldots\}; \quad d=1,2,\ldots;\quad \sum\limits_{i=1}^d(p_i+q_i+1)=n.
$$
 Here all $\alpha_i, \beta_i$ for $i>d$ are zeros.

The measure of such a point equals
$$ 
P^{(n)}_{zz'}(p|q)=
\frac{n!t^d}{(t)_n}\prod\limits
_{i=1}^d{\frac{(z+1)_{p_i}(z'+1)_{p_i}(-z+1)_{q_i}(-z'+1)_{q_i}}{(p_i!)^2(q_i!)^2}}
\cdot {\det}^2\left(\frac{1}{p_i+q_j+1}\right)
$$
where $t=zz'$ and $(a)_k=\Gamma(a+k)/\Gamma(a)$ stands for the Pohgammer symbol. For the parameters $z$ and $z'$ there are two  possibilities, if we assume that the measure of a point described above is always strictly positive,

(1)  $z'=\bar z$ and $z\notin \Bbb Z$

(2)  $z$ and $z'$ are real and $m<z,z'<m+1$ for some $m\in  \Bbb Z$.

One can prove that we really obtain a sequence of probability measures that weakly converges as $n \to \infty$ to some probability measure $P_{zz'}$ on $\Omega$.

The measure $P^{(n)}_{zz'}$ can be regarded as a probability distribution on the set of all
Young diagrams of order $n$, if we consider the numbers $(p_1,\ldots,p_d;q_1,\ldots,q_d)$ as the
Frobenius coordinates of a Young diagram, see [M, \S I.1] for the Frobenius notation.
Furthermore, the values $P^{(n)}_{zz'}(p|q)$ can be regarded as `moments' of the limit measure $P_{zz'}$. Namely,
$$
{P^{(n)}_{zz'}(p|q)}={\dim \lambda}\int_{\Omega}\widetilde{s}_{\lambda}(\omega)P_{zz'}(d\omega)
$$
where $(p|q)$ is considered as the Frobenius notation for a Young diagram $\lambda$; $\dim \lambda$ is the dimension of the complex irreducible representation of the symmetric group $S_n$ corresponding to $\lambda$ (see [M], [JK]); and $\widetilde{s}_{\lambda}$ stands for the so-called {\it extended Schur function}, see [O], [KOO] for details. These moments lie in the base of all our computations.
 
Now we shall introduce the point process, see [DVJ] for general information about point processes. Let us denote by $I$ the punctured interval $[0,1]\setminus \{0\}$. We shall call a set that does not intersect a sufficiently small interval $(-\varepsilon,\varepsilon)$ a {\it test set}. A {\it configuration} is, by definition, a finite or countable system of points in  $I$ such that its intersection with any test set is finite. For any point of $\Omega$ we define the corresponding configuration $(\alpha_1,\alpha_2,\ldots,-\beta_1,-\beta_2,\ldots)$, where we suppose that $\alpha_i$ and $\beta_i$ are nonzero. Then the probability measure $P_{zz'}$ on $\Omega$ provides us a probability measure on the configurations, i.e. a point process. We compute the correlation functions of these processes.

In what follows we use the notion of  the correlation measure as well as that of the correlation function. The latter is defined as the density of the former with respect to the Lebesgue measure. As a distribution it always exists. Note that in [DVJ] the correlation measure is referred to as {\it `factorial moment measure'.} 

We shall work with distributions all over this paper. Everywhere
below the word `distribution' is understood as a synonym of `generalized function', and our notation for the distributions coincides with the usual one for generalized functions.

The paper is organized as follows.

Chapter 1 contains the computation of  the correlation functions. In [O] the computation was reduced to a rather complicated multidimensional moment problem, which we formulate in Section 1.1. This moment problem contains the values of irreducible characters of finite symmetric groups. In Section 1.2 we give a reformulation of the well-known Murnaghan-Nakayama rule for computing these values. The results of Section 1.2 allow to reduce the moment problem to a finite number of simpler moment problems; this reduction is carried out in Section 1.3. In Section 1.4 we introduce three simple operations on compactly supported distributions, these operations correspond to some natural operations on the moments of distributions. Using these operations, in Section 1.5 we construct the solutions of simplified moment problems stated in 1.3. In this construction we use certain distributions $I_d, d=0,1,2,\ldots$ with given moments. To clarify the idea of construction in 1.5, we postpone the explicit computation of $I_d$'s to Section 1.6. After this computation we already have an expression of the correlation functions in the form of a finite sum of some known distributions. However, the summation is taken over a complicated combinatorial set; and we substantially simplify the description of this set in Section 1.7.

It turns out to be much easier to describe the correlation functions (and the whole process), if we restrict ourselves to the behaviour of $\alpha_i$'s (or $\beta_i$'s) only, see the definition of the point processes above. We do this in Chapter 2. Section 2.1 explains that in this case only one summand of the complicated finite sum from Chapter 1 is relevant. In Section 2.2 we produce a simple integral representation of this summand; and it awfully resembles that of multivariate hypergeometric functions. In Section 2.3 we introduce {\it Lauricella hypergeometric functions of type B} and in the next section we express the correlation functions via these Lauricella functions explicitly. These two sections are considered to be optional. They contain rather hard technical work, and results obtained there will be used in also optional Sections 4.2 and 4.3. 
In Section 2.5 we prove that our point processes are simple, i.e., that there are no repetitions of $\alpha_i$'s or $\beta_i$'s with probability $1$.  

Chapter 3 is devoted to a modification of our processes, which we call {\it lifted processes}. The idea is that a simple lifting of measures $P_{zz'}$ from $\Omega_0$ to the bigger space of pairs of convergent series with positive terms substantially simplifies the correlation functions. In Section 3.1 we introduce a general formalism of the lifting and show that its application to well-known  Poisson-Dirichlet processes gives just a Poisson process. In Section 3.2 we demonstrate that the lifting reduces multidimensional integral representations of correlation functions obtained before  to determinants built from a function in $2$ variables with a given integral representation.
In the context of Chapter 2 the situation is even more simpler: the lifted correlation functions are given by determinantal formulas with a kernel expressed via Whittaker functions, we prove this in Section 3.3. We call the kernel the {\it Whittaker kernel}, the explicit formula for it is one of our main results. Determinantal formulas for correlation functions provides a connection of our processes with those arising in random matrices and with certain models of mathematical physics. We hope to clarify this connection  in subsequent papers. 

In Chapter 4 we compute the asymptotics of both lifted and non-lifted correlation function at the origin, which happens to be the same. This computation can be considered as an application of our integral representations. But its main destination is to provide a base for considering so-called {\it tail processes}, which describe the behaviour of particles `infinitely close to zero', see [BO] for further explanations. 
It turns out to be much easier to compute the asymptotics in the lifted case: it suffices to work only with the Whittaker kernel. This is done in Section 4.1. Section 4.2 and 4.3 are optional: they, again, contain many technicalities which do not produce surprising results. However, without them this work would be incomplete. In Section 4.2 we construct asymptotic formulas for Lauricella hypergeometric functions, and in Section 4.3 we use them to compute the asymptotics in the non-lifted case. The computation substantially exploits formulas from Section 2.4.

This work was totally inspired by G.~I.~Olshanski. Besides that he set up the problem, numerous stimulating discussions with him played the crucial role for the whole work. I am happy to express him my deep gratitude.

\head 1. General Structure of  The Correlation Functions  \endhead
\subhead 1.1.  The moment problem \endsubhead
As was shown in [O], the calculation of the correlation functions of our processes can be reduced to the calculation of the family of probability measures on cubes ${[-1,1]}^n,\ n=1,2,\ldots,$ with their moments known. These measures will be referred to as {\it controlling measures} and denoted by $\siz(dx_1,\ldots,dx_n)$;  we shall denote their densities with respect to the Lebesgue measure by $\siz(x_1,\ldots,x_n)$. The densities are considered as distributions (with compact supports). Here $z$ and $z'$ are parameters of the process, see Introduction. The moments of the controlling measures are given by the following explicit formulas, see [O, Proposition 3.3],
$$
\align
&\inte_{{[-1,1]}^n}x_1^{l_1}\cdots x_n^{l_n}\siz (dx_1,\ldots,dx_n)=
\su_{d=1}^n\su_{\Sb p_1>\ldots>p_d\geq 0\\
q_1>\ldots>q_d\geq 0\\ |p|+|q|+d=|l|+n\endSb}
\chi^{(p_1,\ldots,p_d|q_1,\ldots,q_d)}_{(l_1+1,\ldots,l_n+1)}\\
&\times
\frac{t^d}{(t)_{|p|+|q|+d}}\pro_{i=1}^d\frac{(z+1)_{p_i}(-z+1)_{p_i}(z'+1)_{q_i}(-z'+1)_{q_i}}
{p_i!q_i!}
\cdot \det \left(\frac{1}{p_i+q_j+1}\right)
\tag 1.1
\endalign
$$
where $l_1,\ldots,l_n=0,1,2,\ldots$;
$$
|p|=p_1+\ldots+p_d,\quad
|q|=q_1+\ldots+q_d,\quad
|l|=l_1+\ldots+l_n+n;
$$
$(p_1,\ldots,p_d|q_1,\ldots,q_d)$ is the Frobenius notation for a Young diagram, see [M, \S I.1]; and $\chi^{\lambda}_{m_1,\ldots,m_k}$ ($\lambda$ is an arbitrary Young diagram, $m_1,\ldots,m_k$ are positive integers)  stands for the value of the irreducible character  $\chi^{\lambda}$ of the symmetric group $S_{|\lambda|}$ at any permutation with $k$ cycles of length $m_1,\ldots,m_k$, see [JK], [M].

The connection between the correlation and controlling measures is provided by the following
statement proved in  [O, Proposition 4.3].
\proclaim{Theorem 1.1.1}
 On the set $\{(x_1,\ldots,x_n)\in I^n|\ x_i\ne x_j\  \text{for}\  i\ne j\}$
$$
\rz(x_1,\ldots,x_n)=\frac{1}{|x_1\cdots x_n|}\siz(x_1,\ldots,x_n).
$$
\endproclaim
Moreover, using [O, Proposition 4.2] we shall demonstrate (see Section 2.5) that the set $$\{(x_1,\ldots,x_n)\in I^n|\ x_i= x_j\  \text{for some}\  i\ne j\}$$ is a null set with respect to the $n$th correlation measure.

Thus, our problem is completely reduced to the calculation of the controlling measures.
\subhead 1.2. Combinatorics of the Murnaghan-Nakayama rule \endsubhead
In the following we shall need a reformulation of the well-known Murnaghan-Nakayama rule for computing the values of  irreducible characters of  finite symmetric groups.  
\proclaim{Theorem  1.2.1 (Murnaghan-Nakayama rule, see [M, \S I.7, Ex.5])}
Let $\lambda$ be a Young diagram, $|\lambda|=n$; $\rho=(\rho_1,\ldots,\rho_k)$, all $\rho_i$
are positive integers, $\rho_1+\ldots\rho_k=n$. Then 
$$
\chi^{\lambda}_{\rho}=\su_S(-1)^{\operatorname{ht}(S)}
$$
where the sum is taken over all sequences of Young diagrams $S=(\lambda^{(0)}\subset\lambda^{(1)}\subset\ldots\subset\lambda^{(k)}=\lambda)$ such that
$\lambda^{(i)}-\lambda^{(i-1)}$ is a skew $\rho_i$-hook; $\operatorname{ht}(S)$ is the height of $S$,
$$
\operatorname{ht}(S)=\su_i\operatorname{ht}(\lambda^{(i)}-\lambda^{(i-1)}).
$$
(The height of a skew hook is one less than the number of rows it occupies.)
\endproclaim
Let us introduce some notation.

We shall consider certain positive integral-valued variables which will be called {\it linear horizontal} and {\it linear vertical blocks}.

We also use variables whose possible values are pairs $(p,q)\in\{0,1,\ldots\}^2$, and we call these variables {\it hook blocks}.
  
The {\it filling number} of a linear block is its value. For a hook block we shall call $p$ the {\it filling p-number} and $q$ the {\it filling q-number} of this block.

The value of a linear block will be also called the {\it cardinality} of the block. For a hook block we define its cardinality as the sum of its two filling numbers plus one.

A {\it fragment} is, by definition, a partially ordered finite set of blocks subject to the conditions

(1) there is exactly one hook block that precedes all the others;

(2) all linear horizontal blocks are linearly ordered;

(3) all linear vertical blocks are linearly ordered;

(4) there are no other order relations.

Finally, we shall call a finite set of fragments with linear ordering of all their blocks (the ordering is compatible with the partial orderings inside the fragments) a {\it structure}.

A structure is called {\it filled} if all filling numbers of its blocks are known. 

The ordered set of cardinalities of all blocks of a filled structure is called the cardinality of the filled structure (the order of blocks' cardinalities is the same as the order of blocks in the structure) . 

We shall denote by $\upsilon(T)$ the number of linear vertical blocks in $T$. 

Let us define the filling p-number  (q-number) of a fragment as the sum of the p-number of its hook block and the filling numbers of its linear horizontal blocks (the sum of the q-number of the hook block and the filling numbers of its linear vertical blocks, respectively).

Thus, any filled structure $T$ with $d$ fragments produces a set of pairs
$$
\{(P_1,Q_1),\ldots,(P_d,Q_d)\}
$$
 which consists of the filling p- and q-numbers of the fragments. 

Let us define the sign of $T$ as follows
$$
\operatorname{sgn}T=\operatorname{sgn}\Bigl(\pro_{i<j}(P_i-P_j)(Q_i-Q_j)\Bigr)
{(-1)}^{\sum Q_j+\upsilon(T)}.
$$ 
Our main statement in this section is the following
\proclaim{Theorem 1.2.1'}
Under the conditions of Theorem 1.2.1
$$
\chi^{\lambda}_{\rho}=\su_T\operatorname{sgn} T,\tag1.2
$$
where the sum is taken over all filled structures with the cardinality $\rho=(\rho_1,\ldots,\rho_k)$ such that the sets $\{P_1,\ldots,P_d\}$ and $\{Q_1,\ldots,Q_d\}$ of the filling p- and q-numbers of the fragments coincide, up to a permutation, with the Frobenius coordinates of the Young diagram $\lambda$.
\endproclaim
\example{Remark 1.2.2} For the calculation of the $n$th correlation function we shall consider the structures with $n$ blocks. For $n=1$ there is only one such structure; it has one hook block. For $n=2$ we have already 3 such structures: the structure with two fragments, each of them consists of one (hook) block; the structure with one hook and one vertical block; and the structure with one hook and one horizontal block.
\endexample
\demo{Proof of Theorem 1.2.1'} We shall call the structures, which can be obtained from a given one by throwing out several last (in the sense of the ordering inside the structure) blocks,  the {\it substructures} of the given structure. A filling of the initial structure induces the fillings of all its substructures. We shall call a filled structure {\it proper} if for all its substructures (including itself) the following assertion holds: all filling p-numbers (as well as q-numbers) of the fragments are mutually distinct.

The statement of Theorem 1.2.1' immediately follows from the following lemmas.
\proclaim{Lemma 1.2.3}
All proper filled structures $T$ which participate in the sum (1.2) are in one-to-one correspondence with the sequences $S$ from Theorem 1.2.1.  Moreover,
$$
\operatorname{sgn}T={(-1)}^{\operatorname{ht}(S)}\tag 1.3
$$
for corresponding $T$ and $S$.
\endproclaim
\proclaim{Lemma 1.2.4}
The part of the sum (1.2) corresponding to all improper filled structures $T$ vanishes.
\endproclaim 
Clearly, these lemmas and Theorem 1.2.1 imply Theorem 1.2.1'.
\demo{Proof of Lemma 1.2.3} Any removal  of a skew $\theta$-hook from a Young diagram is equivalent either to a removal of two Frobenius coordinates of the diagram with their sum equal to $\theta-1$, or to a reduction of a Frobenius $p$-coordinate or $q$-coordinate by $\theta$. 

More precisely, if we remove a $\theta$-hook $s$ that does not intersect the diagonal and lies to the right of it occupying the rows from $i$th to $j$th, $i<j$, then it is equivalent to the following changes in the sequence of Frobenius $p$-coordinates $(p_1>\ldots>p_d)$ of our diagram
$$
(p_1>\ldots>p_d)\to (p_1>\ldots>p_{i-1}>p_{i+1}>\ldots>p_j>p_i-\theta>p_{j+1}>\ldots>p_d).
$$
By definition, 
$$
(-1)^{\operatorname{ht}(s)}=(-1)^{j-i}.
$$ 

If we remove a $\theta$-hook $s$ that lies completely to the left of the diagonal occupying the columns from $i'$th to $j'$th, then, by analogy, it is equivalent to the following changes in the sequence of Frobenius $q$-coordinates of the diagram
$$
(q_1>\ldots>q_d)\to (q_1>\ldots>q_{i'-1}>q_{i'+1}>\ldots>q_{j'}>q_{i'}-\theta>q_{j'+1}>\ldots>q_d)
$$
The cardinality $\theta$ of $s$ is the sum of the number of columns occupied by $s$ plus the number of rows occupied by $s$ minus one. Therefore, 
$$
(-1)^{\operatorname{ht}(s)}=(-1)^{\theta+j'-i'+1}.
$$ 

Finally,
if we remove a $\theta$-hook that intersects the diagonal of the diagram, occupies rows starting from the $i$th one, and columns starting from the $i'$th one; it is equivalent to removing the $i$th Frobenius $p$-coordinate and the $i'$th Frobenius $q$-coordinate, and 
$$
p_i+q_{i'}=\theta-1.
$$ 
This $\theta$-hook can be divided into three parts: $p_i$-hook lying to the right of the diagonal, $q_{i'}$-hook lying to the left of the diagonal, and a box lying on the diagonal. Using this and previous considerations, one can easily compute that
$$
(-1)^{\operatorname{ht}(s)}=(-1)^{i'-i+q_{i'}}.
$$ 

  We define the correspondence from the hypothesis of Lemma 1.2.3 in such way that three operations described above correspond to deleting a linear horizontal, a linear vertical, or a hook block of a structure respectively. 

Formally it can be made by induction on $k$ (the number of blocks in $T$ or the length of $\rho$). Namely, for $k=1$ the sequence $\emptyset=\lambda^{(0)}\subset \lambda^{(1)}=\lambda$, $\lambda=(p|q)$, corresponds to the structure $T$ with one hook block, and its p- and q-filling numbers are exactly $p$ and $q$. Furthermore, if a sequence $S'$ with $k=m$ is obtained from some sequence $S$ with $k=m+1$ by throwing out the last skew hook $\lambda^{(m+1)}-\lambda^{(m)}$, then the corresponding to $S'$ structure $T'$ is obtained from corresponding to $S$ structure $T$ by throwing out the $(m+1)$th block. If the hook thrown out lies to the right of the diagonal, to the left of the diagonal, or intersects the diagonal, then the corresponding block is linear horizontal, linear vertical, or hook respectively. The Frobenius $p$- and $q$-coordinates of $\lambda=\lambda^{(m+1)}$ coincide (up to permutations) with p- and q-filling numbers of fragments of $T$. 

 Let us check the coincidence of signs (1.3). We shall do this by induction on the number of blocks. If our structure has only one (hook) block with $p$ and $q$ as its filling numbers then the corresponding diagram $\lambda=(p|q)$ is a skew hook itself, and the corresponding sequence $S$ is $\emptyset=\lambda^{(0)}\subset \lambda^{(1)}=\lambda$. In this case 
$$
\operatorname{sgn}T={(-1)}^{\operatorname{ht}(S)}=(-1)^q.
$$
Let us consider a structure $T$ with $m+1$ blocks such that the $(m+1)$th block of $T$ is linear horizontal and its filling number is $\theta$. Let us denote by $T'$ the structure obtained by removing the $(m+1)$th block from $T$ (as above). At the beginning of the proof we considered the corresponding removal of a skew $\theta$-hook $s$. The p-filling numbers $(P_1,\ldots,P_d)$ of fragments of $T$ are permuted numbers $(p_1>\ldots>p_d)$. It is easy to see, that if we reduce one of $P$'s corresponding to $p_i$ by $\theta$, then the sign of $\prod_{i<j}(P_i-P_j)$ and $\operatorname{sgn}T$ will be multiplied by $(-1)^{i-j}$, that  coincides with the sign $(-1)^{\operatorname{ht}(s)}$. So we have reduced the coincidence of signs for $T$ and corresponding sequence $S$ to that for $T'$ and corresponding sequence $S'$, and the number of blocks in $T'$ is less then that in $T$. One can easily check that the same reduction works for the last block of $T$ being linear vertical or hook as well. Modulo this check, (1.3) is proved by induction. 

  The condition that all structures must be proper is necessary because the Frobenius $p$- or $q$-coordinates of any Young diagram are pairwise distinct.\qed
\enddemo   
\demo{Proof of Lemma 1.2.4} Let us exclude from our consideration improper filled structures with equal p- or q-numbers of fragments; their signs are zeros. All the other improper structures will be divided into pairs with opposite signs.

Let us fix one such structure. Let us choose from its substructures the biggest with equal p- or q-numbers. To be concrete, let there be two equal filling p-numbers of some two fragments. Let us exchange in this fragments all linear horizontal blocks which distinguish the initial structure and the substructure. The set of such blocks is not empty, because the initial structure has different p-numbers. Thus we obtained a new filled structure which will be the pair for the initial one.  Obviously, the procedure described above is involute, and the structures in one pair have opposite signs (the Vandermonde determinant  $\pro_{i<j}(P_i-P_j)$ (or  $\pro_{i<j}(Q_i-Q_j)$ in case of q-numbers) changes its sign). \qed
\enddemo
\enddemo
\example{Remark 1.2.5}
The coincidence of  signs (1.3)  for a structure without linear
blocks was communicated to me by G.~I.~Olshanski.
\endexample
\subhead 1.3. A simplification of the moment problem \endsubhead
Our goal in this section is to represent the $n$th correlation measure as a sum of distributions with compact supports;  the summands will be naturally parametrized by (not filled) structures with $n$ blocks. The distributions will be defined as the solutions of  moment problems which are much simpler than the initial one, see 1.1.
\proclaim{Proposition 1.3.1} The moments of the $n$th controlling measure have the form
$$
\multline
\langle\siz(x_1,\ldots,x_n),x_1^{l_1}\cdots x_n^{l_n}\rangle
=\su_T t^d (-1)^{\sum Q_i+\upsilon(T)}\\ \times \frac{\pro_i(z+1)_{P_i}
(z'+1)_{P_i}(-z+1)_{Q_i}(-z'+1)_{Q_i}}
{(t)_{\sum(P_i+Q_i+1)}\pro_i P_i!Q_i!}
\det\left(\frac{1}{P_i+Q_j+1}\right)
\endmultline
$$
where the summation is taken over all filled structures $T$ of the cardinality $(l_1+1,\ldots,l_n+1)$; $(P_1,\ldots,P_d;Q_1,\ldots,Q_d)$ are the filling numbers of the fragments of $T$, $d$ is the number of the fragments.
\endproclaim
\demo{Proof} This statement is an immediate corollary of  the formula (1.1) and Theorem 1.2.1'.
The only thing to mention here is that due to the well known formula for the Cauchy determinant (see [W])
$$
\det\left(\frac{1}{P_i+Q_j+1}\right)=\frac{\pro_{i<j}\bigl((P_i-P_j)(Q_i-Q_j)\bigr)}{\pro_{i,j}(P_i+Q_j+1)}
$$
 we may use the equality
$$
\operatorname{sgn}\det\left(\frac{1}{P_i+Q_j+1}\right)=
\operatorname{sgn}\pro_{i<j}\bigl((P_i-P_j)(Q_i-Q_j)\bigr).\qed
$$
\enddemo
Let us represent the sum from Proposition 1.3.1 as the double sum: the outer summation will be taken over all structures with $n$ blocks, and the inner summation will be taken over all possible fillings of the structure given by the outer sum. Thus we obtain partitions of all the moments into summands parametrized by the unfilled structures with $n$ blocks. In Section 1.5 we shall construct the distribution with a compact support with its moments equal to the summands corresponding to a fixed structure $T$; this distribution will be denoted by $\siz(T)$.  Now we are going to reformulate the moment problem for $\siz(T)$.

Let us fix a structure $T$ with $n$ blocks and $d$ fragments. Let us enumerate the fragments arbitrarily by the numbers from $1$ to $d$. We denote the number of  linear horizontal blocks in the $i$th fragment by $\mu_i$ and the number of linear vertical blocks -- by $\nu_i$.  Thus
$$
\su_{i=1}^d(\mu_i+\nu_i+1)=n
$$
(in each fragment there is exactly one hook block).

In this notation $\upsilon(T)=\sum_{i=1}^d\nu_i$.

We say that the variable $x_j$ corresponds to a block if the block has the number $j$ in the total blocks' ordering of $T$.  Let us rename the variables $x_1,\ldots,x_n$ by the letters 
$$
a_i,\quad \{b_i^j\}_{j=1}^{\mu_i},\quad \{c_i^j\}_{j=1}^{\nu_i};\quad i=1,\ldots,d;
$$
so that the variable $a_i$ corresponds to the hook block of the $i$th fragment;  $\{b_i^j\}_{j=1}^{\mu_i}$ correspond to the horizontal blocks of the $i$th fragment; and $\{c_i^j\}_{j=1}^{\nu_i}$ correspond to the vertical blocks of the $i$th fragment.

By analogy with this, let us rename the numbers $(l_1,\ldots,l_n)$ by the letters $(\{A_i\},\{B_i^j\},\{C^j_i\})$.
\proclaim{Proposition 1.3.2}
With the preceding notation
$$
\siz=\su_T\siz(T),
$$
the sum is taken over all structures $T$ with $n$ blocks and 
$$
\aligned
\langle\siz(T),\pro_{i=1}^d
\Bigl(a_i^{A_i}(b_i^1)^{B_i^1}\cdots(&b_i^{\mu_i})^{B_i^{\mu_i}}(c_i^1)^{C_i^1}\cdots(c_i^{\nu_i})^{C_i^{\nu_i}}\Bigr)
\rangle
\\=t^d\su_{\Sb A_i'+A_i''=A_i\\i=1,\ldots,d\endSb}(-1)^{\sum Q_i+\upsilon(T)}&\cdot\frac{\pro_i(z+1)_{P_i}
(z'+1)_{P_i}
(-z+1)_{Q_i}
(-z'+1)_{Q_i}}
{(t)_{\sum(P_i+Q_i+1)}\pro_iP_i!Q_i!}
\\
&\qquad\qquad\qquad\qquad\times
\det\left(\frac{1}{P_i+Q_j+1}\right)
\endaligned
\tag1.4
$$
where
$$
\cases
P_i=A_i'+\su_{j=1}^{\mu_i}(B_i^j+1)\\
Q_i=A_i''+\su_{j=1}^{\nu_i}(C_i^j+1).
\endcases
\tag1.5
$$
\endproclaim
\demo{Proof}
The only freedom in filling a structure $T$ with a fixed cardinality is that we can vary  filling p- and q-numbers of every hook block with their sum fixed. If we denote the p-number of the hook block in the $i$th fragment of $T$ by $A_i'$ and the q-number of the same block by $A_i''$, then we observe that the cardinality of the block in question equals $A_i+1=A_i'+A_i''+1$.  Computing the filling numbers of the fragments (denoted by $P_i$ and $Q_i$) together with the definition of $\siz(T)$ and Proposition 1.3.1 complete the proof. \qed 
\enddemo
\subhead 1.4. Contraction, pseudoconvolution, and diagonalization \endsubhead
 In this section we shall introduce three operations on distributions which will be used further for the construction of  solutions of our moment problems (see previous section for the setting of the problems).
\proclaim{Proposition 1.4.1 (`Contraction')}
Let $f(\xi_1,\xi_2,\eta_1,\ldots,\eta_k)$, $k\geq 0$, be a distribution with  compact support and the moments
$$
\langle f, {\xi_1}^{\alpha_1}{\xi_2}^{\alpha_2}{\eta_1}^{\beta_1}\ldots{\eta_k}^{\beta_k}\rangle
=m_{\alpha_1\alpha_2\beta_1\ldots\beta_k}.
$$
Then there exists a distribution $g(\xi,\eta_1,\ldots,\eta_k)$ with a compact support and the moments
$$
\langle g, {\xi}^{\alpha}{\eta_1}^{\beta_1}\ldots{\eta_k}^{\beta_k}\rangle
=\su_{\alpha_1+\alpha_2=\alpha}(-1)^{\alpha_2}m_{\alpha_1\alpha_2\beta_1\ldots\beta_k}.
$$

If, moreover, $\operatorname{supp}\subset\{\xi_1\ge 0,\ \xi_2\ge 0\}$, then in the domain $\{\xi\ne 0\}$ the following equality holds
$$
g(\xi,\eta_1,\ldots,\eta_k)=\inte_{y\ge 0}\frac{|\xi|}{|\xi|+y}(f(\xi,y,\eta_1,\ldots,\eta_k)+f(y,-\xi,\eta_1,\ldots,\eta_k))dy.
$$
We shall use the notation 
$$
g(\xi,\eta_1,\ldots,\eta_k)=\Cal C^{\xi}_{\xi_1\xi_2}\left[(f(\xi_1,\xi_2,\eta_1,\ldots,\eta_k)\right].
$$
\endproclaim
\demo{Proof} The value of the distribution $g$ on any test function $\psi(\xi,\eta_1,\ldots,\eta_k)$ is given by the formula
$$
\langle g,\psi\rangle
=\langle f,\frac{\xi_1\psi(\xi_1,\eta_1,\ldots,\eta_k)+ 
\xi_2\psi(-\xi_2,\eta_1,\ldots,\eta_k)}{\xi_1+\xi_2}
\rangle.
$$
\enddemo
\proclaim{Proposition 1.4.2 (`Pseudoconvolution'), see [O]}

Let $f(\xi_1,\ldots,\xi_m)$ and  
$g(\eta_1, \ldots,\eta_m)$ be distributions with compact supports
and the moments
$$
\langle f,{\xi_1}^{\alpha_1}\cdots,{\xi_m}^{\alpha_m}\rangle=
m_{\alpha_1 \ldots\alpha_m},\quad
\langle g,{\eta_1}^{\beta_1}\cdots,{\eta_m}^{\beta_m}\rangle=
n_{\beta_1 \ldots\beta_m}.
$$
Then there exists a distribution $h(\zeta_1,\ldots,\zeta_m)$ with  compact support and the moments
$$
\langle h,{\zeta_1}^{\gamma_1}\cdots,{\zeta_m}^{\gamma_m}\rangle=
m_{\gamma_1 \ldots\gamma_m}n_{\gamma_1 \ldots\gamma_m}.
$$
Moreover,
$$
\operatorname{supp}h\subset \operatorname{supp}f\cdot \operatorname{supp}g.
$$
We shall use the notation 
$$
h=f\odot g.
$$
\endproclaim
\demo{Proof}
The value of the distribution $h$ on a test function $\psi$ is given by the formula
$$
\langle h,\psi\rangle=\int f(\xi_1,\ldots,\xi_m)g(\eta_1,\ldots,\eta_m)\psi(\xi_1\eta_1,\ldots,\xi_m\eta_m)\pro_{i=1}^m d\xi_i d\eta_i.$$
\enddemo
\proclaim{Proposition 1.4.3 (`Diagonalization')}
Let $f(\xi,\eta_1,\ldots,\eta_k)$, $k\geq 0$, be a distribution with the moments
$$
\langle f,\xi^{\alpha}{\eta_1}^{\beta_1}\cdots{\eta_k}^{\beta_k}\rangle=m_{\alpha\beta_1\ldots\beta_k}.
$$
Then the moments of the distribution 
$$
g(\xi_1,\xi_2,\eta_1,\ldots,\eta_k)=f(\xi_1,\eta_1,\ldots,\eta_k)\delta(\xi_2-\xi_1)
$$
equal
$$
\langle g,{\xi_1}^{\alpha_1}{\xi_2}^{\alpha_2}{\eta_1}^{\beta_1}\cdots{\eta_k}^{\beta_k}\rangle=
m_{\alpha_1+\alpha_2\beta_1\ldots\beta_k}.
$$
We shall use the notation 
$$
g(\xi_1,\xi_2,\eta_1,\ldots,\eta_k)=\Cal D^{\xi_1\xi_2}_{\xi}\left[f(\xi,\eta_1,\ldots,\eta_k)\right].
$$
\endproclaim
\demo{Proof} Obvious.
\enddemo
\subhead 1.5.  The distributions $\siz(T)$ \endsubhead
In the present section we shall construct the solutions of the moment problems (1.4) for all $T$ starting from a certain family of distributions. We shall also establish that only few of these solutions really contribute to the correlation functions. 

The crucial role in the whole construction is played by the family of distributions called $I_d(a_1',a_1'';\ldots;a_d,a_d'')$, $d=1,2,\ldots$, whose moments, by definition, are equal to
$$
\align
\langle I_d(a_1',a_1'';&\ldots;a_d',a_d''),(a_1')^{A_1'}(a_1'')^{A_1''}\cdots (a_d')^{A_d'}(a_d'')^{A_d''}\rangle\\=&
\frac{\prod\limits
_i (z+1)_{A_i'}(z'+1)_{A_i'}(-z+1)_{A_i'}(-z'+1)_{A_i''}}{\prod\limits_i A_i'! A_i''! (t)_{\sum{(A_i'+A_i''+1)}}} 
\cdot {\det}\left(\frac{1}{A_i'+A_j''+1}\right).
\tag 1.6
\endalign
$$
The following statement will be proved in the next section, and now, in order to clarify the construction, we shall assume that it is true.
\proclaim{Proposition 1.5.1} There exist distributions $I_d(a_1',a_1'';\ldots;a_d,a_d'')$, $d=1,2,\ldots$, such that their moments are given by (1.6). Moreover,
$$
\operatorname{supp}I_d\subset\{(a_1',a_1'';\ldots;a_d';a_d'')\in \Bbb R_+^{2d}|\su_{i=1}^d(a_i'+a_i'')\leq 1\}.
\tag 1.7
$$
\endproclaim
We use the notation introduced in 1.3. Let us denote
(for all $i=1,\ldots,d$)
$$
\aligned
&\Cal C_i=\Cal C^{a_i}_{a_i'a_i''},\\
 \Cal D_i'=\pro_{j=1}^{\mu_i}\Cal D_{a_i'}^{a_i'b_i^j}&,\quad \Cal D_i''=\pro_{j=1}^{\nu_i}\Cal D_{a_i''}^{a_i'',-c_i^j},\\
&\Cal D_i=\Cal D_i'\Cal D_i''.
\endaligned
\tag1.8
$$
Moreover, let us denote by $\Cal M_i$ (for all $i=1,\ldots,d$) the operation of multiplication of a distribution by $(a_i')^{\mu_i}(a_i'')^{\nu_i}$.
\proclaim{Proposition 1.5.2} For any fixed structure T the compactly supported distribution $\siz(T)$ in variables
$$
a_i,\quad \{b_i^j\}_{j=1}^{\mu_i},\quad \{c_i^j\}_{j=1}^{\nu_i};\quad i=1,\ldots,d;
$$
defined by
$$
\siz(T)=t^d\left(\pro_{i=1}^d\Cal C_i\right)\left(\pro_{i=1}^d\Cal D_i\right)
\left(\pro_{i=1}^d\Cal M_i\right)\cdot I_d (a_1',a_1'';\ldots;a_d',a_d'')
$$
gives a solution of the `simplified moment problem' (1.4).
\endproclaim 
\demo{Proof}
The following proof  is nothing but  sequential applications of statements from the previous section. Let us introduce the distributions
$$
X= \left(\pro_{i=1}^d\Cal M_i\right)\cdot I_d (a_1',a_1'';\ldots;a_d',a_d'')
$$
in variables $a_i',a_i'',\ i=1,\ldots,d$;
$$
Y=\left(\pro_{i=1}^d\Cal D_i\right)\cdot X
$$
in variables $a_i',a_i'',\  \{b_i^j\}_{j=1}^{\mu_i},\ \{c_i^j\}_{j=1}^{\nu_i},\ i=1,\ldots,d;$ and
$$
Z=\left(\pro_{i=1}^d\Cal C_i\right)\cdot Y
$$
in variables  $a_i,\  \{b_i^j\}_{j=1}^{\mu_i},\ \{c_i^j\}_{j=1}^{\nu_i},\ i=1,\ldots,d.$
Let us introduce the following notation. For a distribution $f(\xi_1,\ldots,\xi_k)$ we shall denote 
by $m_{l_1,\ldots,l_k}(f)$ its $(l_1,\ldots,l_k)$-moment. That is,
$$
m_{l_1,\ldots,l_k}(f)=\langle f,\xi_1^{l_1}\cdots \xi_k^{l_k}\rangle.
$$
Denote 
$$
A'=\{A_i'\}_{i=1}^d,\quad A''=\{A_i''\}_{i=1}^d,\quad A=\{A_i\}_{i=1}^d;
$$
$$
B=\{B_i^j\}_{i=1,j=1}^{d,\mu_i},\quad C=\{C_i^j\}_{i=1,j=1}^{d,\nu_i};
$$
$$
\mu=\{\mu_i\}_{i=1}^d,\quad \nu=\{\nu_i\}_{i=1}^d;
$$
cf. 1.3. Then, by Proposition 1.4.1,
$$
m_{A,B,C}(Z)=\su_{\Sb A'+A''=A\\
A_i',A_i''\geq 0\endSb}(-1)^{|A''|}m_{A',A'',B,C}(Y)
\tag 1.9
$$
where the absolute value sign for a vector stands for the sum of  its coordinates.  
Furthermore, by Proposition 4.3 we get
$$
m_{A',A'',B,C}(Y)=(-1)^{|\overline{C}|}m_{A'+\overline{B},A''+\overline{C}}(X)
\tag 1.10
$$
where
$$
\overline{B}=\{\overline {B_i} \}_{i=1}^d,\quad \overline{B_i}=\su_{j=1}^{\mu_i}B_i^j;
$$
$$
\overline{C}=\{\overline{C_i}\}_{i=1}^d,\quad \overline{C_i}=\su_{j=1}^{\nu_i}C_i^j, \quad |\overline{C}|=\su_{i,j}C_i^j.
$$
Note that the sign $(-1)^{|\overline{C}|}$ appears by the following fact: for a one dimensional distribution $g(\eta)$
$$
\langle g(\eta),\eta^l\rangle =(-1)^{l}\langle g(-\eta),\eta^l\rangle.
$$
 In our case we apply this fact to all variables $\{c_i^j\}$, see the minus sign before $c_i^j$ in (1.8).

Finally,
$$
m_{A'+\overline{B},A''+\overline{C}}(X)=m_{A'+\overline{B}+\mu,A''+\overline{C}+\nu}(I_d).
\tag 1.11
$$
Recall the notation (1.5). It can be reformulated as
$$
\align
&P=A'+\overline{B}+\mu,\\
&Q=A''+\overline{C}+\nu.
\endalign
$$ 
In particular, $(-1)^{|A''|+|\overline{C}|}=(-1)^{|Q|+|\nu|}=(-1)^{|Q|+\upsilon(T)}$. Combining (1.9), (1.10), (1.11) together
we obtain (using (1.6))
$$
\multline
m_{A,B,C}(Z)=\su_{\Sb A'+A''=A\\ A_i',A_i''\geq 0\endSb}(-1)^{|Q|+\upsilon(T)}m_{P,Q}(I_d)=\su_{\Sb A'+A''=A\\ A_i',A_i''\geq 0\endSb}\frac{(-1)^{|Q|}}{(t)_{|P|+|Q|+d}}\\
\times\prod\limits_{i=1}^d\frac{ (z+1)_{P_i}(z'+1)_{P_i}(-z+1)_{Q_i}(-z'+1)_{Q_i}}{ P_i! Q_i! } 
\cdot {\det}\left(\frac{1}{P_i+Q_j+1}\right).
\endmultline
$$
Hence, the moments of $t^dZ$ coincide with those from (1.4), as was to be proved. \qed
\enddemo
Let us call a structure {\it substantial} if there is at most one linear block in each its fragment. In the preceding notation it means that $\mu_i+\nu_i\leq1$ for all $i$. The next statement shows (see also Theorem 1.1.1 and a few words after it) that only substantial structures are relevant.
\proclaim{Proposition 1.5.3} If a structure $T$ is not substantial then
$$
\operatorname{supp}\siz(T)\subset \{(x_1,\ldots,x_n)|\pro_i x_i\cdot\pro_{i<j}(x_i-x_j)=0\}.
$$
\endproclaim
\proclaim{Lemma 1.5.4} Let the $(\alpha_1,\ldots,\alpha_k)$-moment of a distribution  $f(\xi_1,\ldots,\xi_k)$ with  compact support depend on $\alpha_1,\alpha_2$ only via their sum.
Then
$$
\operatorname{supp}f\subset\{(\xi_1,\ldots,\xi_n)|\ \xi_1=\xi_2\}.
$$
\endproclaim
\demo{Proof of Lemma 1.5.4} $(\xi_1-\xi_2)f(\xi_1,\ldots,\xi_k)=0$.\enddemo
\demo{Proof of Proposition 1.5.3} Let us use the preceding notation. If a structure $T$ is not substantial then there are three possible cases: $\mu_i\geq 2$ for some $i$,  $\nu_i\geq 2$ for some $i$, or there exists such $i$ that $\mu_i=\nu_i=1$.

 Let us first assume that $d=1$, i.e, our structure consists of only one fragment. We shall omit the subscript $i$ of all variables and operations introduced above, because in this case $i$ is identically equal to one. 

Suppose $\mu\geq 2$. Then the moments of $\siz(T)$, see (1.4), depend on $B^1$ and $B^2$ only via their sum (because they depend on $P=B^1+B^2+\ldots$).  Thus, Lemma 1.5.4 implies the assertion of the theorem.

Suppose $\nu\geq 2$. Then, again, the moments of $\siz(T)$, see (1.4), depend on $C^1$ and $C^2$ only via their sum (because they depend on $Q=C^1+C^2+\ldots$). Again Lemma 1.5.4 completes the proof.

Finally, suppose $\mu=\nu=1$. By Proposition 1.5.2 we know that
$$
\Cal M I_1(a',a'')=X,\quad \Cal D X=\Cal D'D'' X=Y,\quad \Cal CY=Z,
$$
$$ 
\siz(T)=tZ.
$$
We have 
$$
X=\Cal M I_1(a',a'')=a'a''I_1(a',a''),\quad \Cal D''X=a'a''I_1(a',a'')\delta(c+a''),
$$
$$
Y=\Cal D'D'' X=a'a''I_1(a',a'')\delta(c+a'')\delta(b-a').
$$
Furthermore, applying Proposition 1.4.1 (we can use the explicit formula for the contraction because the support of $I_1(a',a'')$ lies in the domain $a',a''\geq 0$, see Proposition 1.5.1), for $a\neq 0$ we obtain
$$
\multline
Z(a,b,c)=\Cal C\cdot  Y(a',a'',b,c)=\int_{y\geq 0}\frac{|a|}{|a|+y}(Y(a,y,b,c)+Y(y,-a,b,c))dy\\
=\int_{y\geq 0}\frac{|a|}{|a|+y}\bigl(ayI_1(a,y)\delta(c+y)\delta(b-a)-ayI_1(y,-a)\delta(c-a)\delta(b-y)\bigr)dy\\
=-\frac{a^2c}{a-c}I_1(a,-c)\delta(b-a)+\frac{a^2b}{b-a}I_1(b,-a)\delta(c-a).
\endmultline
$$
Obviously, the support of the last expression lies in the domain 
$$
\{(a,b,c)|\ (b-a)(c-a)=0\} ,
$$
as was to be proved. 

In case of arbitrary $d$ we use the fact that the operations $\Cal C_i,\Cal D_i,\Cal M_i$ with different subscripts act on different sets of variables. That is why the same considerations as above applied to a fragment with linear vertical and linear horizontal blocks  prove the assertion. \qed 
\enddemo

Now let us fix a numeration of the fragments of a substantial structure $T$ (this numeration  was arbitrary up to this moment) in such a way that
$$
\align
&\mu_i=1,\ \nu_i=0;\quad i=1,\ldots,m_1,\\
&\mu_i=0,\ \nu_i=1;\quad i=m_1+1,\ldots,m_1+m_2,\\ 
&\mu_i=0,\ \nu_i=0;\quad i=m_1+m_2+1,\ldots,m_1+m_2+m_3=d.
\endalign
$$
In other words, we denoted the number of the fragments with horizontal blocks by $m_1$ and put them first, then we put $m_2$ fragments with vertical blocks, and, finally, we have $m_3$ fragments without linear blocks at all. 

According to Theorem 1.1.1 we  also set
$$
\rho_n^{(zz')}(T)=\frac{1}{|x_1\cdots x_n|}\siz(T).
\tag1.12 
$$

It turns out that the case $m_3=0$ is much simpler than the general situation. Namely, in this case the formulas for $\siz(T)$ and for $\rz(T)$ are rather simple. 
\proclaim{Theorem 1.5.5} Let $m_3=0$. Then in the domain where all the variables are nonzero and pairwise distinct 
$$
\rho_{2d}^{(zz')}(T)=t^d\pro_{i=1}^n\frac{1}{r_i+|s_i|} I_d(r_1,-s_1;\ldots,r_d,-s_d)
\tag 1.13
$$
where
$$
\align
&(r_1,\ldots,r_d)=(b_1,\ldots,b_{m_1},a_{m_1+1},\ldots,a_d);\\
&(s_1,\ldots,s_d)=(a_1,\ldots,a_{m_1},c_{m_1+1},\ldots,c_d).
\endalign
$$
\endproclaim
\example {Remark 1.5.6} Theorem 1.5.5 and Proposition 1.5.1 imply that for $m_3=0$ $\rz(T)$ is a distribution defined in the domain where the sums of the absolute values of variables corresponding to the two blocks of each fragment are nonzero (each fragment contains exactly two blocks).  Moreover,  Theorem 1.5.5 implies that the support of $\rho_n^{(zz')}(T)$ lies (note that we have thrown out the set $\prod_{i<j}(x_i-x_j)=0$) in the set where all variables corresponding to linear horizontal blocks (i.e. $\{b_i\}$) are nonnegative, all variables corresponding to linear vertical blocks (i.e. $\{c_i\}$) are nonpositive, and every two variables corresponding to the two blocks of the same fragment (i.e. $a_i$ and $b_i$ for $i\leq m_1$, or $a_i$ and $c_i$ for $m_1<i\leq m_1+m_2$) have different signs.
\endexample 
\demo{Proof of Theorem 1.5.5}
 We shall work out a detailed proof in two cases: $(m_1=1,\ m_2=0)$, and $(m_1=0,\ m_2=1)$, and then give an explanation how to get the proof in the general case.

Suppose $m_1=1$ and $m_2=0$. We shall follow all the steps demonstrated in the proof of Proposition 1.5.2. By this proof we know that (cf. the proof of Proposition 1.5.3)
$$
\Cal M I_1(a',a'')=X,\quad \Cal D X=\Cal D' X=Y,\quad \Cal CY=Z,
$$
$$ 
\sigma^{(zz')}_2(T)=tZ.
$$
We have 
$$
X=\Cal M I_1(a',a'')=a'I_1(a',a''),\quad Y=\Cal D'X=a'I_1(a',a'')\delta(b-a').
$$
Using Proposition 1.4.1, namely, the explicit formula for the contraction, for $a\neq 0$ we get
$$
\multline
Z(a,b)=\Cal C_{a'a''}^a\cdot  Y(a',a'',b)=\int_{y\geq 0}\frac{|a|}{|a|+y}(Y(a,y,b)+Y(y,-a,b))dy\\
=\int_{y\geq 0}\frac{|a|}{|a|+y}\bigl(aI_1(a,y)\delta(b-a)+yI_1(y,-a)\delta(b-y)\bigr)dy\\
=\int_{y\geq 0}\frac{a^2}{a+y}I_1(a,y)dy\cdot\delta(b-a)-\frac{ab}{b-a}I_1(b,-a).
\endmultline
$$
Note that the support of the first summand lies in the set where $a=b$. Thus, we have showed that  in the domain where $a\neq 0$ and $a\neq b$
$$
\rho_2^{(zz')}(T)=tZ(a,b)=t \frac{|a|b}{b+|a|}I_1(b,-a),
$$
as was to be proved (in this case $r_1=b$ and $s_1=a$). 

Suppose $m_1=0$ and $m_2=1$. Then
$$
\Cal M I_1(a',a'')=X,\quad \Cal D X=\Cal D'' X=Y,\quad \Cal CY=Z,\quad
\sigma^{(zz')}_2(T)=tZ;
$$
and
$$
X=\Cal M I_1(a',a'')=a''I_1(a',a''),\quad Y=\Cal D''X=a''I_1(a',a'')\delta(c+a'').
$$
Again, by Proposition 1.4.1, when $a\neq 0$ we get
$$
\multline
Z(a,c)=\Cal C_{a'a''}^a\cdot  Y(a',a'',c)=\int_{y\geq 0}\frac{|a|}{|a|+y}(Y(a,y,c)+Y(y,-a,c))dy\\
=\int_{y\geq 0}\frac{|a|}{|a|+y}\bigl(yI_1(a,y)\delta(c+y)-aI_1(y,-a)\delta(c-a)\bigr)dy\\
=-\frac{ac}{a-c}I_1(a,-c)+\int_{y\geq 0}\frac{a^2}{y-a}I_1(y,-a)dy\cdot\delta(c-a).
\endmultline
$$
Now the support of the second summand lies in the set where $c=a$, and hence in the domain where $a\neq 0$ and $a\neq c$ we have
$$
\rho_2^{(zz')}(T)=tZ(a,c)=t \frac{a|c|}{a+|c|}I_1(a,-c).
$$
In this case $r_1=a$ and $s_1=c$.

In the general case for  arbitrary number of fragments $d$ we can apply the considerations demonstrated above to each of the fragments. Indeed, the operators $\Cal C_i,\ \Cal D_i,\ \Cal M_i$ act on different variables if they have different subscripts. Moreover, each of the fragments is of one of the types considered above: it either contains a linear horizontal block (like in the case $m_1=1,\ m_2=0$) or a linear vertical block  (like in the case $m_1=0,\ m_2=1$). Thus, in general case the proof is obtained by word for word applications of one of our two previous considerations to the appropriate fragments: the first part applies to the first $m_1$ fragments, and the second one -- to the last $m_2$ fragments. \qed  
\enddemo
The case $m_3>0$ is a bit more complicated. As we have seen in the proofs of Proposition 1.5.3 and Theorem 1.5.5, the considerations for different fragments are independent. That is why, in order to understand the general situation, let us consider the unique structure that consists of only one (hook) block. That is, in our notation, $m_1=m_2=0,\ m_3=1$. Then
$$
\Cal M I_1(a',a'')=I(a',a'')=X,\quad \Cal D X=X=I_1(a',a'')=Y,
$$
$$
\Cal CY=Z,\quad
\sigma^{(zz')}_1(T)=tZ.
$$
By Proposition 1.4.1, for $a\neq 0$
$$
\multline
Z(a)=\Cal C_{a'a''}^a\cdot  Y(a',a'')=\int_{y\geq 0}\frac{|a|}{|a|+y}(Y(a,y)+Y(y,-a))dy\\
=\int_{y\geq 0}\frac{|a|}{|a|+y}\bigl(I_1(a,y)+I_1(y,-a)\bigr)dy,
\endmultline
$$
and in the same domain
$$
\rho_1^{(zz')}(a)=t\int_{y\geq 0}\frac{I_1(a,y)}{a+y}dy+t\int_{y\geq 0}\frac{I_1(y,-a)}{|a|+y}dy.
\tag 1.14
$$
(Note that we have completely calculated the first correlation function, because there exists only one structure with one block.) Proposition 1.5.1 implies that the supports of the two summands do not overlap: in the first summand $a>0$, while in the second  $a<0$. 

The answer (1.14) can be reformulated in the following way. Consider all substantial structures  with 1 fragment, which contain the initial structure as a substructure, and such that for them $m_3=0$, i.e.,  they have a linear (vertical or horizontal) block. In our case there are exactly two such structures -- those that were considered in the proof of Theorem 1.5.5. Let us fix one of these two structures and denote it by $\widehat T$. Then it contains  one extra block with respect to the initial structure $T$. We take the distribution  $\rho_2^{(zz')}(\widehat T)$ described by Theorem 1.5.5 and integrate it over the variable corresponding to the extra linear block. Let us denote the answer by $\rho_1^{(zz')}(T,\widehat T)$. Then  we have proved (formula (1.14)) that 
$$
\rho_1^{(zz')}(T)=\su_{\widehat T} \rho_1^{(zz')}(T,\widehat T)
\tag 1.15
$$
where the sum is taken over all `enveloping' structure described above (i.e., the sum contains two summands). As was mentioned before, the supports of the two summands do not overlap.

This complicated explanation of (1.14) has only one advantage: an analog of formula (1.15) holds in general case. 

Let $T$ be an arbitrary substantial structure. We fix the numeration of its fragments as described before Theorem 1.5.5. 

We shall call $\widehat{T}$ an {\it enveloping} structure of $T$ if $\widehat{T}$ is substantial; $T$ is a substructure of $\widehat{T}$; $T$ and $\widehat{T}$ have the same number of fragments; the numeration of the extra blocks in $\widehat{T}$ (with respect to $T$) is compatible with the numeration of fragments of $T$; in each fragment of $\widehat{T}$ there is exactly one linear block.  

In other words, in order to obtain an enveloping structure of $T$ we have to add in each fragment of $T$ with only one (hook) block a linear (horizontal or vertical) block, observing the ordering of fragments. The number of enveloping structures equals, obviously, $2^{m_3}$
where, as before, $m_3$ is the number of fragments of $T$ without linear blocks.  

For $\widehat{T}$ being an enveloping structure of $T$ we define
$$
\rho_n^{(zz')}(T,\widehat{T})=\inte_{x_{n+1},\ldots,x_{n+m_3}}\rho_{n+m_3}^{(zz')}(\widehat{T})
dx_{n+1}\cdots dx_{n+m_3}.
\tag 1.16
$$

I.e., $\rho_n^{(zz')}(T,\widehat{T})$ is obtained from $\rho_{n+m_3}^{(zz')}(\widehat{T})$ by integrating the latter over all variables which correspond to the extra $m_3$ linear blocks added to $T$ for obtaining $\widehat{T}$, and $\rho_{n+m_3}^{(zz')}(\widehat{T})$ is given by the general formula (1.13) (by definition, any enveloping structure has no fragments with only one block).
\example{Remark 1.5.7} The distribution $\rho_n^{(zz')}(T,\widehat{T})$ does not depend on the numeration of fragments in $T$ fixed before. A change of this numeration is equivalent to a change of the numeration of extra $m_3$ variables $(x_{n+1},\ldots,x_{n+m_3})$ that does not affect the result of the integration in (1.16).  
\endexample
\proclaim{Theorem 1.5.8} In the domain where all variables are nonzero and pairwise distinct
$$
\rho_n^{(zz')}(T)=\su_{\widehat{T}}\rho_n^{(zz')}(T,\widehat{T})
\tag 1.17
$$
where the sum is taken over all enveloping structures $\widehat{T}$ of $T$. 
\endproclaim
\example{Remark 1.5.9} Using Remark 1.5.6 and (1.16) it is easy to see that the supports of  different summands of (1.17) do not overlap (in the domain where all variables are nonzero and mutually distinct), the support of each summand lies in one of the hyperoctants $\{x_i > \text{or} <0,i=1,\ldots,n\}$, and the hyperoctants are different for different summands. Thus, the support of $\rho_n^{(zz')}(T)$ lies in the union of $2^{m_3}$ hyperoctants which are described by the following conditions (cf. Remark 1.5.6): all variables corresponding to linear horizontal blocks  are positive, all variables corresponding to linear vertical blocks are negative, and every two variables corresponding to  two blocks of the same fragment  have different signs. 
\endexample
\demo{Proof of Theorem 1.5.8} We have already done all the work. One have to apply the considerations, which led to the formulas (1.14) and (1.15), to each of $m_3$ fragments without linear blocks. The other $m_1+m_2$ fragments are considered as in Theorem 1.5.5. As we know, different fragments are treated independently, and the result of all these applications will be exactly the formula (1.17).\qed
\enddemo
 In Theorems 1.5.8 and 1.5.5 we have completely determined (modulo the construction of $I_d$, see next section) the correlation functions in the domain where all variables are nonzero and mutually distinct, using the notion of structures. But the combinatorics of substantial structures and their enveloping structures seems rather complicated. In Section 1.7 we shall interpret the pairs $(T,\widehat{T})$ as the elements of a simpler combinatorial object.
\subhead 1.6. Construction of $I_d$\endsubhead
The goal of the present section is to prove Proposition 1.5.1 and to provide explicit formulas for the distributions $I_d,\ d=1,2,\ldots$. We shall start using this formulas in Chapter 2.

Let us introduce the distributions 
$$
\phi_a(u)=\frac{u^a_+}{\Gamma (a+1)},\quad u\in \Bbb R,\ a\in \Bbb C.
$$
 For $\Re a>-1$,  $\phi_a$, by definition, equals $u^a/\Gamma(a+1)$ for $u>0$ and vanishes for $u<0$, so it is a locally integrable function. For $\Re a\leq -1$, $\phi_a$ is defined via analytic continuation. For example, $\phi_{-1}(u)=\delta(u)$.

We shall also deal with products of the type
$$
\phi_{ab}=\phi_a(u)\phi_b(1-u),\quad a,b\in \Bbb C,
$$
which are well defined for the reason that possible singularities of the factors do not overlap.

Let us also remind that in Proposition 1.4.2 we have defined an operation $\odot$ on compactly supported distributions in the same number of variables and called it {\it pseudoconvolution}. 
The characteristic property of the pseudoconvolution of two distributions is that its moments equal the products of moments of initial distributions.
\proclaim{Proposition 1.6.1} For every $d=1,2,\ldots$ there exists a compactly supported distribution $J_d(a_1',a_1'';\ldots;a_d',a_d'')$ with moments
$$
\align
\langle J_d(a_1',a_1'';&\ldots;a_d',a_d''),(a_1')^{A_1'}(a_1'')^{A_1''}\cdots (a_d')^{A_d'}(a_d'')^{A_d''}\rangle\\=&
\frac{\prod\limits
_i (z+1)_{A_i'}(z'+1)_{A_i'}(-z+1)_{A_i'}(-z'+1)_{A_i''}}{\prod\limits_i A_i'! A_i''! (t)_{\sum{(A_i'+A_i''+1)}}} 
\pro_{i=1}^d\frac{1}{A_i'+A_j''+1}.
\tag 1.18
\endalign
$$
Moreover,
$$
\operatorname{supp}J_d\subset\{(a_1',a_1'';\ldots;a_d';a_d'')\in \Bbb R_+^{2d}|\su_{i=1}^d(a_i'+a_i'')\leq 1\}.
\tag 1.19
$$
\endproclaim
\demo{Proof}
Using pseudoconvolution, let us define the distributions $J_d$ as follows
$$
J_d= F_1\odot F_2\odot F_3,
\tag 1.20
$$
where
$$
F_1(a_1',a_1'';\ldots;a_d',a_d'')=\Gamma(t)\prod\limits_{i=1}^d \phi_{z'}(a_i')\phi_{-z'}(a_i'')\cdot \phi_{t-d-1}\bigl(1-\sum\limits_{i=1}^d (a_i'+a_i'')\bigr);
$$
$$
F_2(a_1',a_1'';\ldots;a_d',a_d'')=\prod\limits_{i=1}^d \phi_{z,-z-1}(a_i')\phi_{-z,z-1}(a_i'');
$$
$$
F_3(a_1',a_1'';\ldots;a_d',a_d'')=\delta (a_1'-a_1'',\ldots, a_d'-a_d'') \prod\limits_{i=1}^d \chi_{[0,1]}(a_i')
$$
Here
$$
\chi_{[0,1]}(u)=
\cases
1,&\quad  u\in [0,1] \cr
0,&\quad  u\notin [0,1] 
\endcases
$$
is the characteristic function of the segment [0,1].
By the fact that the support of pseudoconvolution is a subset of the pointwise product of the supports of the factors, the inclusion (1.19) follows from the following obvious relations
$$
\operatorname{supp}F_1\subset\{(a_1',a_1'';\ldots;a_d';a_d'')\in \Bbb R_+^{2d}|\su_{i=1}^d(a_i'+a_i'')\leq 1\}.
$$
$$
\operatorname{supp}F_2,\ \operatorname{supp}F_3\subset\{(a_1',a_1'';\ldots;a_d';a_d'')|\ 0\leq a_i',a_i''\leq 1,\ i=1,\ldots,d\}.
$$

Let us compute the moments of $J_d$ and compare them with (1.18). 
We shall use the following lemma which is a well-known generalization of the Euler's beta-integral (the integral below is also called Dirichlet integral). 
\proclaim {Lemma 1.6.2} For any $\alpha_o,\alpha_1,\ldots,\alpha_m\in \Bbb C$ 
$$
\inte_{\Sb u_1,\ldots,u_n\\ u_0=1-u_1-\ldots-u_m\endSb}\phi_{\alpha_0}(u_0)\cdots\phi_{\alpha_m}(u_m)du_1\cdots du_m=\frac{1}{\Gamma(\alpha_0+\ldots+\alpha_m+m)}.
$$
\endproclaim
\demo{Proof of Lemma 1.6.2} Induction on $m$. \enddemo

Note that by this statement we can immediately compute the moments of $F_1$, $F_2$, and $F_3$. Indeed, the moments of $F_1$ is $2d$--dimensional Dirichlet integral. We get
$$
\langle F_1, \pro_{i=1}^d\bigl((a_i')^{A_i'}(a_i'')^{A_i''}\bigr)\rangle=
\frac{\pro_{i=1}^d (z'+1)_{A_i'}(-z'+1)_{A_i''}}{(t)_{\sum(A_i'+A_i''+1)}}.
$$
The moments of $F_2$ are products of  $2d$ Euler beta-integrals,
$$
\langle F_2, \pro_{i=1}^d\bigl((a_i')^{A_i'}(a_i'')^{A_i''}\bigr)\rangle=
\pro_{i=1}^d \frac{(z+1)_{A_i'}(-z+1)_{A_i''}}{A_i'!A_i''!}.
$$
Finally, the moments of $F_3$ are the products of $d$ integrals of monomials $(a_i')^{A_i'+A_i''}$ over $[0,1]$. Thus
$$
\langle F_3, \pro_{i=1}^d\bigl((a_i')^{A_i'}(a_i'')^{A_i''}\bigr)\rangle=
\pro_{i=1}^d\frac{1}{A_i'+A_i''+1}.
$$
The moments of pseudoconvolution are products of moments of factors, see Proposition 1.4.2. Thus, we obtain (1.18).\qed
\enddemo
\demo{Proof of Proposition 1.5.1}
We set
$$
I_d(a_1',a_1'';\ldots;a_d',a_d'')=\su_{\sigma\in S_d}\operatorname{sgn}\sigma\cdot
J_d(a_1',a_{\sigma(1)}'';\ldots;a_d',a_{\sigma(d)}'').
\tag 1.21
$$
Then (1.18) implies (1.6), and (1.19) implies (1.7).\qed
\enddemo
By the fact that $F_1$ and $F_2$ are invariant under the permutations of $\{a_i''\}$, we can write down the following formula for $I_d$, which easily follows from (1.20), (1.21).
\proclaim{Corollary 1.6.3} For every $d=1,2,\ldots$
$$
I_d=\sum\limits_{\sigma\in S_d} \operatorname{sgn} \sigma\cdot F_1\odot F_2\odot F_3^{\sigma},
\tag 1.22
$$
where $F_1,F_2$ are as above, and
$$
\aligned
F_3^{\sigma}(a_1',a_1'';\ldots;a_d',a_d'')=&\delta (a_1'-a_{\sigma(1)}'',\ldots, a_d'-a_{\sigma(d)}'') \prod\limits_{i=1}^d \chi_{[0,1]}(a_i')\\
=&F_3(a_1',a_{\sigma^{-1}(1)}'';\ldots;a_d',a_{\sigma^{-1}(d)}''). 
\endaligned
$$
Here
$$
\chi_{[0,1]}(u)=
\cases
1,&\quad  u\in [0,1] \cr
0,&\quad  u\notin [0,1] 
\endcases
$$
is the characteristic function of the segment [0,1].
\endproclaim
\example{Example 1.6.4} By the formulas above,
$$
I_1(a',a'')=J_1(a',a'')=f_1\odot f_2\odot f_3,
$$
where
$$
\align
f_1(a',a'')&=\Gamma(t) \phi_{z'}(a')\phi_{-z'}(a'')\cdot \phi_{t-2}\bigl(1- (a'+a'')\bigr);
\\
f_2(a',a'')&= \phi_{z,-z-1}(a')\phi_{-z,z-1}(a'');
\\
f_3(a',a'')&=\delta (a'-a'')  \chi_{[0,1]}(a').
\endalign
$$
This statement can be rewritten in the following form.
$$
\align
I_1(a',a'')=\Gamma(t)\inte_{u,v} &\phi_{z'}(u)\phi_{-z'}(v)\cdot \phi_{t-2}(1- u-v)\\ \times& \inte_{w\in [0,1]}\phi_{z,-z-1}\left(\frac{a'}{uw}\right)\phi_{-z,z-1}\left(\frac{a''}{vw}\right)
\frac{dudvdw}{uvw^2}.
\tag 1.23
\endalign
$$
Indeed, the value of such distribution on a test function $\psi(a',a'')$ equals, by definition,
$$
\align
\Gamma(t)\inte_{u,v} \phi_{z'}(u)\phi_{-z'}(v)\cdot \phi_{t-2}(1- u-v)\qquad\qquad\qquad\qquad\qquad\qquad\qquad&\\ \times \inte_{w\in [0,1]}\left[\inte_{a',a''}\phi_{z,-z-1}(a')\phi_{-z,z-1}(a'')\psi(a'uw,a''vw)da'da''\right]&
dudvdw\\=\langle f_1\odot f_2&\odot f_3,\psi\rangle,
\tag 1.24
\endalign
$$
see 1.4.2. 

If we want to interpret the integrals in (1.23) over $u,v,w$ as the usual ones and not as a formal sign of pairing, we need to impose some conditions; in particular,  the set $\{(u,v,w)|\ uvw=0\}$ has to be negligible. From (1.24) we see that if $f_1$ is an integrable function, then everything works (the expression in brackets is a smooth function in $u,v,w$). Thus, in order to understand the integration in (1.23) as the usual one, we may require, for example, $-1<\Re z'<1,\ t>2$. We shall discuss these problems in more details at the beginning of Chapter 2, where we shall use the multidimensional analog of (1.23).  
\endexample
\subhead 1.7. Combinatorics of the substantial structures \endsubhead
In this section we give a nice combinatorial reformulation of Theorem 1.5.8.

Let us denote by $\Phi_{n,d}$ the set of mappings
$$
\varphi :\{1,\ldots,n\}\to \{1,1';\ldots;d,d'\}
$$
subject to the two conditions

1) $\varphi$ is injective, i.e. $\varphi(i)\ne \varphi(j)$ if $i\ne j$;

2) $\operatorname{Im} \varphi \cap\{m,m'\}\ne\emptyset$ for all $m=1,\ldots,d.$

It is clear that $\Phi_{n,d} \ne \emptyset$ if and only if $n/2\le d\le n$. On the set $\{1,1';\ldots;d,d'\}$ we have a natural action of the symmetric group $S_d$: this group permutes the pairs $(i,i')$. The action induces an action of $S_d$ on $\Phi_{n,d}$, and the condition 2) implies that every orbit of this action consists  of exactly $|S_d|=d!$ points.
\proclaim{Proposition 1.7.1} There exists a map $\gamma_{n,d}$ of the set $\Phi_{n,d}$ onto the set of all pairs $(T,\widehat{T})$ where $T$ is a substantial structure with $n$ blocks and $d$ fragments, and $\widehat{T}$ is an enveloping structure of $T$. The inverse image of any pair $(T,\widehat{T})$ with respect to this map consists of exactly one orbit of $S_d$ in $\Phi_{n,d}$.
\endproclaim
\demo{Proof}
Let us construct $\gamma_n,d$. We fix $\varphi\in \Phi_{n,d}$ and produce the pair $(T,\widehat{T})=\gamma_{n,d}(\varphi)$ as follows. Let us describe the $i$th fragment of $T$ and $\widehat{T}$. 

(1) If $\varphi^{-1}(i)=\emptyset$ then the $i$th fragment of $T$ has only one (necessarily hook) block, and the number of its block equals $\varphi^{-1}(i')$. The enveloping structure $\widehat{T}$ has in this fragment one extra linear horizontal block.

(2) If $\varphi^{-1}(i')=\emptyset$ then the $i$th fragment of $T$ also has only one (necessarily hook) block, and the number of its block equals $\varphi^{-1}(i)$. The enveloping structure $\widehat{T}$ has in this fragment one extra linear vertical block.

(3) If $\varphi^{-1}(i)>\varphi^{-1}(i')$ then the $i$th fragment of $T$ has one hook block number $\varphi^{-1}(i')$ and one linear horizontal block number $\varphi^{-1}(i)$. The $i$th fragment of $\widehat{T}$ coincides with that of $T$.

(4) If $\varphi^{-1}(i')>\varphi^{-1}(i)$ then the $i$th fragment of $T$ has one hook block number $\varphi^{-1}(i)$ and one linear vertical block number $\varphi^{-1}(i')$. The $i$th fragment of $\widehat{T}$ coincides with that of $T$.

Clearly, $\gamma_{n,d}$ is a surjection. Moreover, given a substantial structure $T$ (with $n$ blocks and $d$ fragments) together with a numeration of its fragments and its enveloping structure $\widehat{T}$, we can restore $\varphi$ using the conditions 1)-4) uniquely. A change of  the numeration of $d$ fragments of $T$ exactly corresponds to the action of $S_d$ on $\Phi_{n,d}$ described above. This numeration is not determined by $T$ or $\widehat{T}$, thus, the inverse image of any pair $(T,\widehat{T})$ is exactly one orbit of $S_d$. \qed
\enddemo

 Let  us introduce some notation. Starting from a function (or a distribution) $F(r_1,s_1;\ldots;r_d,s_d)$ in $2d$ variables and a map $\varphi\in\Phi_{n,d}$ we define the function $(\varphi F)(x_1,\ldots, x_n)$ in $n$ variables, $\varphi \in \Phi_{n,d}$, as follows.
Let us rename the variable $r_i$ by $x_k$ if $\varphi (k)=i$, $s_j$ by $x_k$ if $\varphi (k)=j'$, and let us do this for all $r_i,s_j$ such that $i,j'\in \operatorname{Im} \varphi$. Then let us integrate $F$ over all $r_l,s_m$ such that $l \notin \operatorname{Im} \varphi, \ m'\notin \operatorname{Im} \varphi$. The result is denoted by $(\varphi F)(x_1,\ldots, x_n)$.

For the sake of convenience we also introduce distributions ($d=1,2,\ldots$)
$$
H_d(r_1,s_1;\ldots;r_d,s_d)=\frac{t^d}{d!\prod\limits_{i=1}^d(r_i+|s_i|)}I_d(r_1,-s_1;\ldots;r_d,-s_d)
\tag 1.25
$$
where $I_d$ is defined by (1.22), cf (1.13).
These distributions are defined in the domain where $\prod\limits_{i=1}^d(r_i+|s_i|)\neq 0$. Proposition 1.5.1 implies that 
$$
\operatorname{supp}H_d\subset\{(r_1,s_1;\ldots;r_d;s_d)\in \Bbb R^{2d}|r_i\geq 0,\ s_i\leq 0,\ \su_{i=1}^d(r_i+|s_i|)\leq 1\}.
\tag 1.26
$$
 
The distributions $I_d(r_1,s_1;\ldots;r_d,s_d)$  are symmetric under the permutations of  pairs $(r_i,s_i)$ (because their moments (1.6) are symmetric). Therefore, by (1.25), $H_d(r_1,s_1;\ldots;r_d,s_d)$ are also symmetric under these permutations.

\proclaim{Proposition 1.7.2} For any $\varphi\in\Phi_{n,d}$
$$
(\varphi H_d)(x_1,\ldots,x_n)=\frac{1}{d!}\rho_n^{(zz')}(T,\widehat{T})
$$
where
$$
 (T,\widehat{T})=\gamma_{n,d}(\varphi).
$$
\endproclaim
\demo{Proof}
 Follows from the definitions and the fact that $H_d(r_1,s_1;\ldots;r_d,s_d)$ is symmetric under the permutations of  pairs $(r_i,s_i)$. 
\enddemo
Finally, combining Theorem 1.5.8 and Propositions 1.7.1, 1.7.2, we obtain the main statement of the first chapter.
\proclaim{Theorem 1.7.3} The correlation functions $\rho _n^{(zz')}(x_1,\ldots,x_n)$ in the domain where all variables are nonzero and pairwise distinct have the form
$$
\rho^{(zz')}_n(x_1\ldots,x_n)=\sum\limits_{d\geq n/2}^n \sum \limits_{\varphi \in \Phi_{n,d}}(\varphi H_d)(x_1,\ldots,x_n).
$$
\endproclaim
\head 2. The Correlation Functions in Positive (Negative) Hyperoctants \endhead
By the words 'positive (negative) hyperoctant' we mean the domain where all variables are positive 
(negative). We have the fundamental relation (see [O, Proposition 4.6])
$$
\rho_n^{(zz')}(x_1,\ldots,x_n)=\rho_n^{(-z,-z')}(-x_1,\ldots,-x_n)
\tag2.0
$$
which implies that we may  consider only positive hyperoctants.
\subhead 2.1. Why positive hyperoctants are simpler \endsubhead
Remark 1.5.9 shows that the support of $\rz(T)$ has common points with the positive hyperoctant  only if all $n$ blocks of the substantial structure $T$ are hook blocks. For each $n$ such structure $T$ is unique. Moreover, the only summand of (1.17) which gives a nonzero contribution in the positive hyperoctant is such that the enveloping structure $\widehat{T}$ has only extra (with respect to $T$) linear vertical blocks (thus, no extra horizontal blocks).  There is exactly one such enveloping structure. Therefore, using (1.13), (1.16), and the definition of $H_d$ (1.25) we obtain 
\proclaim{Proposition 2.1.1} For $x_1,\ldots,x_n>0$ in the domain where all variables are pairwise distinct 
$$
\align
\rz(x_1,\ldots,x_n)&=n!\inte_{s_1,\ldots,s_n}H_n(x_1,s_1;\ldots;x_n,s_n)ds_1\cdots ds_n\\
&=t^n\inte_{s_1,\ldots,s_n}\frac{1}{\prod_{i=1}^n(x_i+s_i)}I_n(x_1,s_1;\ldots;x_n,s_n)ds_1\cdots ds_n
\tag 2.1
\endalign
$$
where the distributions $I_n$ were defined in Proposition 1.5.1.
\endproclaim
Thus, we reduced the complicated combinatorial expression for the correlation functions obtained before (Theorem 1.7.3) to the integral (2.1), provided that all our variables are positive.
\subhead 2.2 An integral representation \endsubhead 
In this section we shall derive from Proposition 2.1.1 the following formula.
\proclaim{Theorem 2.2.1} Let $x_1,\ldots,x_n>0$ and $\sum_{i=1}^nx_i<1$. Then in the domain where $\prod_{i<j}(x_i-x_j)\neq 0$,
$$
\align
\rz(x_1,\ldots,x_n)=t^n&\Gamma(t)\inte_{\Sb a_i,b_i\\ i=1,\ldots,n\endSb}\pro_{i=1}^n
\phi_{-z}(a_i)\phi_{z'}(a_i+1)\phi_{-z'}(b_i)\phi_{z}(b_i+1)\\
\times \det &\left(\frac{1}{a_i+b_j+1}\right)\phi_{t-n-1}\bigl(1-\su_{i=1}^nx_i(a_i+b_i+1)\bigr)
\pro_{i=1}^n da_idb_i.
\tag 2.2
\endalign
$$
\endproclaim

\example{Remark 2.2.2}  In Section 2.5 we shall prove that the set where $\prod_{i<j}(x_i-x_j)=0$ is a null set for the correlation functions and that the  restriction $\prod_{i<j}(x_i-x_j)\neq 0$ on the domain in Theorem 2.2.1 is irrelevant.
\endexample
\demo{Proof} 
Note, first of all, that the right-hand side of (2.2) is a well-defined generalized function in our domain. Indeed, the product
$$
\pro_{i=1}^n
\phi_{-z}(a_i)\phi_{z'}(a_i+1)\phi_{-z'}(b_i)\phi_{z}(b_i+1)\cdot \phi_{t-n-1}\bigl(1-\su_{i=1}^nx_i(a_i+b_i+1)\bigr)
$$
is a well-defined distribution in $a_i$'s and $b_i$'s, because any point $(a,b)$ belongs to at most $2n$ singular hyperplanes corresponding to the singularities of factors. The integral in (2.2) is the value of this distribution on the test function
$$
 \det \left(\frac{1}{a_i+b_j+1}\right). 
$$

Secondly, the right-hand side of (2.2) depends on $z$ and $z'$ analytically, and so does the left-hand side (the correlation function). That is why in this proof we may restrict ourselves to the domain  where $-1<\Re z,\Re z'<0$ and $t=zz'>n+1$. Then all integrals appearing in this section can be considered as the usual ones, and not as formal signs of pairing (cf. the end of 1.6). Moreover, all distributions $\phi_a$ in this proof become just usual integrable functions, and we shall operate with them accordingly to this.  
 
By the explicit formulas (1.20), (1.21)  for $I_d$ and the definition of the pseudoconvolution (section 1.4), we obtain (cf. Example 1.6.4)
$$
\align
I_n(r_1,s_1;\ldots;r_n,s_n&)=\Gamma(t)
\inte_{\Sb u_i,v_i\\  i=1,\ldots,n \endSb} \pro_{i=1}^n
\phi_{z'}(u_i)\phi_{-z'}(v_i)\cdot \phi_{t-n-1}\bigl(1-\su_{i=1}^n(u_i+v_i)\bigr)
\\
\times \su_{\sigma\in S_n}&\operatorname{sgn}(\sigma)
 \inte_{\Sb w_i\in [0,1]\\i=1,\ldots,n\endSb} \pro_{i=1}^n\phi_{z,-z-1}\left(\frac{r_i}{u_iw_i}\right)\phi_{-z,z-1}\left(\frac{s_i}{v_iw_{\sigma(i)}}\right)\frac{du_idv_i}{u_iv_iw_i^2}.
\tag 2.3
\endalign
$$ 

Let $f$ be a distribution with compact support, and $\operatorname{supp}f\subset \Bbb R_+$.  
A function $\goth Sf$ of one positive variable is called the  {\it Stieltjes transform} of
 $f$  (cf. [W]), if
$$
\goth Sf(y)=\inte_x \frac{f(x)}{x+y}dx.
$$
Clearly, 
$$
(\goth Sf)(\alpha y)= \goth S (f(\alpha y)),\quad \alpha>0.
$$
Substituting (2.3) into (2.1) and using this observation, we obtain
$$ 
\align
\rz(x_1,&\ldots,x_n)=t^n\Gamma(t)
\inte_{\Sb u_i,v_i\\  i=1,\ldots,n \endSb} \pro_{i=1}^n
\phi_{z'}(u_i)\phi_{-z'}(v_i)\cdot \phi_{t-n-1}\bigl(1-\su_{i=1}^n(u_i +v_i)\bigr)
\\
\times \su_{\sigma\in S_n}&\operatorname{sgn}(\sigma)
 \inte_{\Sb w_i\in [0,1]\\i=1,\ldots,n\endSb} \pro_{i=1}^n\phi_{z,-z-1}\left(\frac{x_i}{u_iw_i}\right)(\goth S\phi_{-z,z-1})\left(\frac{x_i}{v_iw_{\sigma(i)}}\right)\frac{du_idv_idw_i}{u_iv_iw_i^2}.
\tag 2.4
\endalign
$$
We shall need the following 
\proclaim{Lemma 2.2.3}
$$
\goth S \phi_{-z,z-1}(y)=y^{-z}(1+y)^{z-1}.
$$
\endproclaim
\demo{Proof of Lemma 2.2.3} We shall denote by $F(a,b;c;x)$ the standard Gauss hypergeometric function,
$$
F(a,b;c;x)=\Gamma(c)\inte_{\tau}\frac{\phi_{b-1}(t)\phi_{c-b-1}(1-\tau)}{(1-x\tau)^a}d\tau,
$$
 see [E]. We have
$$
\multline
\inte_x\frac{1}{x+y}\phi_{-z}(x)\phi_{z-1}(1-x)dx=\frac{1}{y}\inte_x\frac{\phi_{-z}(x)\phi_{z-1}(1-x)}{1+x/y}dx\\
=\frac{1}{y}F(-z+1,1;1;-1/y)=\frac{1}{y}(1+1/y)^{-z-1}=y^{-z}(1+y)^{z-1}
\endmultline
$$ 
where we used well-known following identities
$$
F(a,b;c;x)=F(b,a;c;x);\qquad F(a,b,b;x)=(1-x)^{-a}.\qed
\tag2.5
$$
\enddemo

By applying the lemma to (2.4) we get
$$
\align
\rz (x_1,\ldots,x_n)=t^n\Gamma (t)\inte_{\Sb u_i,v_i\\  i=1,\ldots,n \endSb} \pro_{i=1}^n
\phi_{z'}(u_i)\phi_{-z'}(v_i)\cdot \phi_{t-n-1}\bigl(1&-\su_{i=1}^n(u_i +v_i)\bigr)\\ \times
 \su_{\sigma\in S_n}\operatorname{sgn}(\sigma)
\inte_{\Sb 0\leq w_i\leq 1\\ i=1,\ldots,n\endSb} \pro_{i=1}^n \phi_{z,-z-1}\left(\frac{x_i}{u_iw_i}\right)\left(\frac{x_i}{v_iw_{\sigma(i)}}\right)^{-z}
&\Bigl(1+\frac{x_i}{v_iw_{\sigma(i)}}\Bigr)^{z-1} \tag 2.6\\ &\times
\frac{du_idv_idw_i}{u_iv_iw_i^2}.
\endalign
$$
We proceed by the following
\proclaim{Lemma 2.2.4}
$$
\inte_0^1 \phi_{z,-z-1}\left(\frac{r_1}{uw}\right)\left(\frac{r_2}{vw}\right)^{-z}
\left(1+\frac{r_2}{vw}\right)^{z-1}
\frac{dw}{uvw^2}=\frac{\phi_{-z}(u-r_1)\phi_z(r_2+v)}{r_1v+r_2u}.
$$
\endproclaim
\demo{Proof of Lemma 2.2.4}
First, let us simplify the expression
$$
\multline
\phi_{z,-z-1}\left(\frac{r_1}{uw}\right)\left(\frac{r_2}{vw}\right)^{-z}
\left(1+\frac{r_2}{vw}\right)^{z-1}
\frac{1}{uvw^2}\\
=\frac{1}{\Gamma(z+1)}\left(\frac{r_1}{uw}\right)^z_+\frac{1}{\Gamma(-z)}\left(1-\frac{r_1}{uw}\right)^{-z-1}_+\left(\frac{r_2}{vw}\right)^{-z}\left(1+\frac{r_2}{vw}\right)^{z-1}\frac{1}{uvw^2}\\
=\frac{(uw-r_1)_+^{-z-1}}{\Gamma(-z)}\frac{(vw+r_2)^{z-1}}{\Gamma(z+1)}.
\endmultline
$$
Let us make the change of variables $w\to y$
$$
y=\frac{uw-r_1}{u-r_1};\quad w=\frac{(u-r_1)y+r_1}{u};\quad dw=\frac{u-r_1}{u}dy.
$$
Then we obtain
$$
\multline
\inte_0^1\frac{(uw-r_1)^{-z-1}_+}{\Gamma(-z)}\frac{(vw+r_2)^{z-1}}{\Gamma(z+1)}dw\\=\inte_0^1\frac{(y(u-r_1))^{-z-1}_+}{\Gamma(-z)}\frac{((y(u-r_1)+r_1)v/u+r_2)^{z-1}}{\Gamma(z+1)}\frac{u-r_1}{u}dy\\
=u^{-z}(u-r_1)^{-z}_+\frac{(vr_1+ur_2)^{z-1}}{\Gamma(z+1)}\inte_0^1\frac{y^{-z-1}}{\Gamma(-z)}\left(1+\frac{v(u-r_1)}{vr_1+ur_2}y\right)^{z-1}dy\\
=u^{-z}(u-r_1)^{-z}_+\frac{(vr_1+ur_2)^{z-1}}{\Gamma(z+1)}\frac{1}{\Gamma(-z+1)}F\left(1-z,-z;1-z;-\frac{v(u-r_1)}{vr_1+ur_2}\right)
\endmultline
$$
by the definition of the hypergeometric function, see above. Using the identities (2.5) we finally get, that our initial integral equals
$$
 u^{-z}\frac{(u-r_1)^{-z}_+}{\Gamma(-z+1)}\frac{(vr_1+ur_2)^{z-1}}{\Gamma(z+1)}
\left(1+\frac{v(u-r_1)}{vr_1+ur_2}\right)^z 
=\frac{(u-r_1)^{-z}_+}{\Gamma(-z+1)}\frac{(v+r_2)^z}{\Gamma(z+1)}\frac{1}{vr_1+ur_2}
$$
The proof of Lemma 2.2.4 is complete.
\enddemo
To complete the proof of Theorem 2.2.1 it only remains to use the result of the previous lemma in (2.6) $n$ times taking
$$
r_1=x_i,\quad r_2=x_{\sigma^{-1}(i)},\quad u=u_i,\quad v=v_{\sigma^{-1}(i)},\quad w=w_i;\quad
i=1,\ldots,n.
$$
Then we have 
$$
\multline
\rz (x_1,\ldots,x_n)=t^n\Gamma (t)\inte_{\Sb u_i,v_i\\  i=1,\ldots,n \endSb} \pro_{i=1}^n
\phi_{z'}(u_i)\phi_{-z'}(v_i)\cdot \phi_{t-n-1}\bigl(1-\su_{i=1}^n(u_i +v_i)\bigr)\\ \times
\su_{\sigma\in S_n}\operatorname{sgn}(\sigma) \pro_{i=1}^n
\phi_{-z}(u_i-x_i)\phi_{z}(v_i+x_i)\frac{1}{x_{\sigma(i)}v_i+x_iu_{\sigma(i)}}du_idv_i
\endmultline
$$
We arrive at the formula (2.2) via the change of variables
$$
\cases
a_i=u_i/x_i-1 \\
b_i=v_i/x_i
\endcases.
\qed
$$
\enddemo
We conclude this section by the following statement which will be used later.
\proclaim{Proposition 2.2.5}
The right-hand side of the formula (2.2) is an analytic function of $x_1,\ldots,x_n$ in the domain 
$$
\{x_1,\ldots,x_n>0;\quad \su_{i=1}^nx_i<1\}.
$$
\endproclaim
\demo{Proof} Set
$$
|x|=x_1+\ldots+x_n.
$$
 By the change of variables
$$
\cases
A_i=\frac{x_i}{1-|x|}a_i\\
B_i=\frac{x_i}{1-|x|}b_i
\endcases
$$
we obtain the following formula for the right-hand side of (2.2)
$$
\multline
t^n\Gamma(t)\frac{(1-|x|)^{t+n-1}}{\prod_{i=1}^nx_i^2}\inte_{\Sb A_i,B_i\\ i=1,\ldots,n\endSb}\pro_{i=1}^n
\phi_{-z}(A_i)\phi_{-z'}(B_i)\phi_{t-n-1}\bigl(1-\su_{i=1}^n(A_i+B_i)\bigr)
\\
\times \phi_{z'}\left(A_i+\frac{x_i}{1-|x|}\right)\phi_{z}\left(B_i+\frac{x_i}{1-|x|}\right)\\
\times \det \left(\frac{1}{\frac{1-|x|}{x_i}A_i+\frac{1-|x|}{x_j}B_j+1}\right)\pro_{i=1}^n dA_idB_i.
\endmultline
$$
This formula can be interpreted as the value of the distribution
$$
\phi_{-z}(A_i)\phi_{-z'}(B_i)\phi_{t-n-1}\bigl(1-\su_{i=1}^n(A_i+B_i)\bigr)
$$
on the test function
$$
\multline
t^n\Gamma(t)\frac{(1-|x|)^{t+n-1}}{\prod_{i=1}^nx_i^2}\phi_{z'}\left(A_i+\frac{x_i}{1-|x|}\right)\phi_{z}\left(B_i+\frac{x_i}{1-|x|}\right)\\
\times
\det \left(\frac{1}{\frac{1-|x|}{x_i}A_i+\frac{1-|x|}{x_j}B_j+1}\right).
\endmultline
$$
The test function is an analytic function of the variables $A_i,B_i$ and of the parameters $x_1,\ldots,x_n$. Clearly, it implies that the value of a distribution that does not depend on the parameters on such test function is also an analytic function in $x_1,\ldots,x_n$. \qed
\enddemo
\example{Remark 2.2.6} Another proof of Proposition 2.2.5 can be obtained from Theorem 2.4.1.
\endexample
Proposition 2.2.5 and Theorem 2.2.1 immediately imply
\proclaim{Corollary 2.2.7} The correlation functions $\rho_n^{(zz')}(x_1,\ldots,x_n)$ are analytic in the domain where all variables are of the same sign, pairwise distinct, and the sum of their absolute values is less then one.
\endproclaim
 In 2.5 we shall see that the restriction $\prod_{i<j}(x_i-x_j)\neq 0$ can be removed (cf. Remark 2.2.2).
\subhead 2.3*. Lauricella hypergeometric functions of type B \endsubhead Let us recall that sections marked by * are considered to be optional, see Introduction for details. 

 One can notice that the integral representation (2.2) of the correlation functions looks like those of multivariate hypergeometric functions. Essentially, the only thing that distinguishes the integral (2.2) from a hypergeometric integral (2.8) is the determinant
$$
\det\left(\frac{1}{a_i+b_j+1}\right).
$$
Our goal in Sections 2.3 and 2.4 is to get rid of this determinant and to give an explicit expression of the correlation functions via so-called Lauricella hypergeometric functions of type B.
Our reference for multivariate hypergeometric functions is the book [AK].

Let $$a=(a_1,\ldots,a_m)\in \Bbb C^m,$$ $$b=(b_1,\ldots,b_m)\in  \Bbb C^m,$$ $$c\in \Bbb C\setminus \{0,-1,\ldots\},$$ $$y=(y_1,\ldots,y_m).$$
Then the $m$-dimensional Lauricella hypergeometric function of type B is defined by the series $$
F_B^{[m]}(a,b;c|y)=\sum\limits_{k_1,\ldots,k_m\geq 0} \frac{(a_1)_{k_1} (b_1)_{k_1}\cdots
(a_m)_{k_m}(b_m)_{k_m}}{(c)_{k_1+\ldots+k_m}k_1!\cdots k_m!}y_1^{k_1}\cdots y_m^{k_m}
\tag 2.7
$$
that is absolutely convergent for $|y_1|<1,\ldots,|y_m|<1$, see [AK].

This function can be analytically continued into the domain $\Re y_1<0,\ldots,\Re y_n<0$ by the Euler-Laplace type integral
$$
F_B^{[m]}(a,b;c|y)=\Gamma(c)\int\limits_{u_1,\ldots,u_m} {\prod\limits_{i}\frac{\phi _{b_i-1}(u_i)}{(1-u_iy_i)^{a_i}} \phi_{c-\sum\limits_{i}b_i-1}\bigl(1-\sum\limits_i u_i\bigr)\prod\limits_i du_i.}
\tag 2.8
$$
 
In the next chapter we shall use another type of integral representations for these functions, namely, Mellin-Barnes type integrals.
\example{Remark 2.3.1}
The function $F_B^{[m]}$ for $m=1$ coincides with the Gauss hypergeometric function (cf. [E]), and for $m=2$ it coincides with the Appell hypergeometric function $F_3$, see [AK], [E].
\endexample
\proclaim{Proposition 2.3.2}
The function  $F_B^{[m]}(a,b;c|y)$ is invariant under the permutations of  triples $(a_i,b_i,y_i)$ 
and under the transpositions $a_i\leftrightarrow b_i$ for all
 $i=1,\ldots,n$.
\endproclaim
\demo{Proof} In the domain $|y_1|<1,\ldots,|y_m|<1$ the claim is obvious from the series (2.7), and the invariance survives after the analytic continuation.\qed
\enddemo
\proclaim{Proposition 2.3.3} For all $k=1,\ldots,m$
$$
\frac{\partial}{\partial y_k}F_B^{[m]}(a,b;c|y)=\frac{a_kb_k}{c}F_B^{[m]}(\hat a,\hat b;c+1|y)
$$
where
$$
\align
\hat a&=(a_1,\ldots,a_{k-1},a_k+1,a_{k+1},\ldots,a_m);\\
\hat b&=(b_1,\ldots,b_{k-1},b_k+1,b_{k+1},\ldots,b_m).
\endalign
$$  
\endproclaim
\demo{Proof} As in the previous proposition, the statement is obvious in the domain where the series (2.7) converges; and the property in question is stable under the analytic continuation.
\qed
\enddemo
The integral representation (2.2) of the $n$th correlation function is given by a 
$2n$-dimensional integral. Hence, it is natural to express $\rho_n^{(zz')}$ via the 
$2n$-dimensional Lauricella function $F_B^{[2n]}$. It turns out to be convenient 
to introduce certain function in $2n$ variables such that 
 $\rho_n^{(zz')}$ will is its restriction to a $n$-dimensional plane, 
 and to express this new function in terms of $F_B^{[2n]}$. We shall do this in 2.4, 
 and now we introduce another function in $2n$ variables, which will be extensively
  used in the next section. 
 
Consider the function $f_n^{(zz')}(y_1',\ldots,y_n';y_1'',\ldots,y_n'')$
in $2n$ variables defined by the following formula 
$$
\align
f_n^{(zz')}(y_1',&\ldots,y_n';y_1'',\ldots,y_n'')\\&=\frac{1}{\prod\limits_i(y_i'-y_i'')}\sum\limits\Sb \epsilon_i=0,1\\ i=1,\ldots,n\endSb (-1)^{\sum\epsilon_i}\prod (y_i')^{1-\epsilon_i} (y_i'')^{\epsilon_i}
 F_B^{[2n]}(a^\epsilon,b^\epsilon;c|y)
\tag2.9
\endalign
$$
where $y_i'\neq y_i''$ for all $i=1,\ldots,n$; and
$$
\aligned
y&=(y_1',\ldots,y_n';y_1'',\ldots,y_n'')\\
a^\epsilon&=(1-\epsilon_1-z',\ldots,1-\epsilon_n-z';\epsilon_1-z,\ldots,\epsilon_n-z)\in \Bbb C^{2n},\\
b^\epsilon&=(1-\epsilon_1-z,\ldots,1-\epsilon_n-z;\epsilon_1-z',\ldots,\epsilon_n-z')\in \Bbb C^{2n},\\
c&=t-n(z+z'-1).
\endaligned
\tag 2.9'
$$
\proclaim{Proposition 2.3.4} The function $f_n^{(zz')}(y_1',\ldots,y_n';y_1'',\ldots,y_n'')$ can be uniquely continued to the points where $y_i'=y_i''$ for some $i=1,\ldots,n$, so that the result will be an analytic function in $y_1',\ldots,y_n';y_1'',\ldots,y_n''$. 
\endproclaim
 The extended function will be also denoted by $f_n^{(zz')}$.
 \demo{Proof}
It suffices to prove that the sum over all $\epsilon_i$ in the right-hand side of (2.9) is skew-symmetric with respect to the transpositions $y_i'\leftrightarrow y_i''$ for all $i=1,\ldots,n$.
For the function $F_B^{[2n]}(a^\epsilon,b^\epsilon;c|y)$, as follows from  Proposition 2.3.2, the transposition $y_i'\leftrightarrow y_i''$ is equivalent to the transposition of the pairs
$$
(a^\epsilon_i,b^\epsilon_i)=(1-\epsilon_i-z',1-\epsilon_i-z)\longleftrightarrow(\epsilon_i-z,\epsilon_i-z')= (a^\epsilon_{n+i},b^\epsilon_{n+i})
$$
that is equivalent to two transpositions $a^\epsilon_i\leftrightarrow b^\epsilon_i$; $a^\epsilon_{n+i}\leftrightarrow b^\epsilon_{n+i}$ (which do not change anything, see Proposition 2.3.2) and the change $\epsilon_i\to 1-\epsilon_i$. For the factor 
$(y_i')^{1-\epsilon_i} (y_i'')^{\epsilon_i}$ the transposition $y_i'\leftrightarrow y_i''$ is also equivalent to the change $\epsilon_i\to 1-\epsilon_i$. But under this change the factor $ (-1)^{\sum\epsilon_i}$ changes its sign.\qed
\enddemo
\example{Example 2.3.5} Let us consider the case $n=1$. Then by applying the L'H\^{o}pital rule to the definition of  $f_n^{(zz')}$ and using Proposition 2.3.3 we obtain
$$
\multline
f_1^{(zz')}(y',y')=F_B^{[2]}(a,b;c|y)\\+y'\Bigl(F_B^{[2]}(a',b';c+1|y)-\frac{zz'}{(1-z)(1-z')}F_B^{[2]} (a'',b'';c+1|y)\Bigr)
\endmultline
$$
where $y=(y',y')$, $c=(1-z)(1-z')$, and
$$
\aligned
a=(1-z',-z),&\qquad b=(1-z,-z');
\\
a'=(2-z',-z),&\qquad b'=(2-z,-z');
\\
a''=(1-z',1-z),&\qquad b''=(1-z,1-z').
\endaligned
\tag2.10
$$
\endexample

\subhead 2.4*. Expression of the correlation functions via $F_B^{[m]}(a,b;c|y)$\endsubhead The main formulas of this section will play a crucial role in Chapter 4.

Now we are introducing the function in $2n$ variables such that the $n$th correlation function is its restriction to a $n$-dimensional plane, see the discussion after Proposition 2.3.3.
   
Let us define the function $R_n^{(zz')}( y_1',\ldots,y_n';y_1'',\ldots,y_n'')$ in $2n$ positive variables by the formula  
$$
\align
R_n^{(zz')}(x_1',&x_1'';\ldots;x_n',x_n'')=t^n\Gamma(t)\inte_{\Sb a_i,b_i\\ i=1,\ldots,n\endSb}\pro_{i=1}^n
\phi_{-z}(a_i)\phi_{z'}(a_i+1)\phi_{-z'}(b_i)\phi_{z}(b_i+1)\\
\times \det &\left(\frac{1}{a_i+b_j+1}\right)\phi_{t-n-1}\Bigl(1-\su_{i=1}^n\bigl(x_i'(a_i+1/2)+x_i''(b_i+1/2)\bigr)\Bigr)
\pro_{i=1}^n da_idb_i.
\tag 2.11
\endalign
$$
Comparing (2.11) and (2.2) we see that if $\prod_{i<j}(x_i-x_j)\neq 0$ then
$$
R_n^{(zz')}(x_1,\ldots,x_n;x_1,\ldots,x_n)=\rz(x_1,\ldots,x_n).
\tag 2.12
$$
Set $c=t-n(z+z'-1)$.
\proclaim{Theorem 2.4.1}  For $x_1,\ldots,x_n>0$ such that $\sum_{i=1}^nx_i<1$ and
 $\prod_{i<j}(x_i-x_j)\neq 0$ 
$$
\align 
\rz(x_1,\ldots,x_n)=\Gamma(t)\prod\limits_{i=1}^n&\phi_{z-1}(x_i)\phi_{z'-1}(x_i)\cdot
\phi_{c-1}(1-|x|)\\
\times &\sum\limits_{\sigma\in S_n}\operatorname{sgn} \sigma \cdot f_n^{(zz')}
(y_1,\ldots,y_n;y_{\sigma(1)},\ldots,y_{\sigma(n)}).
\tag 2.13
\endalign
$$
where $|x|=\sum_{i=1}^nx_i$ and for all $i=1,\ldots,n$
$$
y_i=-\frac{1-|x|}{x_i}.
\tag2.13'
$$
The function $f_n^{(zz')}$ is explicitly expressed via Lauricella hypergeometric functions
in (2.9), (2.9').
\endproclaim
About the restriction $\prod_{i<j}(x_i-x_j)\neq 0$ see Remark 2.2.2.
\demo{Proof} In fact, we shall prove a more general fact. Namely,
if $\sum_{i=1}^n(x_i'+x_i'')<2$ then
$$
\align 
R_n^{(zz')}(x_1',\ldots,x_n';x_1'',\ldots,x_n'')=\Gamma(t)\prod\limits_{i=1}^n&\phi_{z-1}(x'_i)\phi_{z'-1}(x''_i)\cdot\phi_{c-1}\Bigl(1-\frac{ |x'|+|x''|}{2}\Bigr)\\
\times \sum\limits_{\sigma\in S_n}&\operatorname{sgn} \sigma \cdot f_n^{(zz')}(y_1',\ldots,y'_n;y_{\sigma(1)}'',\ldots,y_{\sigma(n)}'')
\tag 2.14
\endalign
$$
where $|x'|=\sum_{i=1}^nx_i'$, $|x''|=\sum_{i=1}^nx_i''$, and for all $i=1,\ldots,n$
$$
y_i'=-\frac{1-\frac{ |x'|+|x''|}{2}}{x_i'},\quad y_i''=-\frac{1-\frac{ |x'|+|x''|}{2}}{x_i''}.
\tag 2.14'
$$
By (2.12), (2.13) is an obvious corollary of (2.14).
 
Observe that (2.11) implies 
$$
\align
R_n^{(zz')}(x_1',\ldots,x_n';x_1'',\ldots,x_n'')&\\=\su_{\sigma\in S_n}\operatorname{sgn}\sigma
\cdot & \widetilde{R}_n^{(zz')}(x_1',\ldots,x_n';x_{\sigma(1)}'',\ldots,x_{\sigma(n)}'').
\tag2.15
\endalign
$$
where
$$
\aligned
\widetilde{R}_n^{(zz')}(x_1'&,\ldots,x_n';x_1'',\ldots,x_n'')\\=t^n&\Gamma(t)\inte_{\Sb a_i,b_i\\ i=1,\ldots,n\endSb}\pro_{i=1}^n
\phi_{-z}(a_i)\phi_{z'}(a_i+1)\phi_{-z'}(b_i)\phi_{z}(b_i+1)\\
\times  &\phi_{t-n-1}\Bigl(1-\su_{i=1}^n\bigl(x_i'(a_i+1/2)+x_i''(b_i+1/2)\bigr)\Bigr)
\pro_{i=1}^n \frac{da_idb_i}{a_i+b_i+1}.
\endaligned
\tag 2.16
$$

(To obtain this it suffices to develop the determinant in (2.11) and renumber the integration variables $\{b_i\}$.)

Next we shall get rid of the unpleasant factor $\prod_i(a_i+b_i+1)^{-1}$ (cf. the beginning of 2.3). 

Note that (we use the obvious formula $\phi_c'=\phi_{c-1}$)
$$
\align
\phi_{t-n-1}\Bigl(1-\su_{i=1}^n\bigl(x_i'(a_i+1/2)+x_i''(b_i+1/2)\bigr)\Bigr)&=
\pro_{i=1}^n\frac{1}{x_i'-x_i''}\\ \times\pro_{i=1}^n\left(\frac{\partial}{\partial b_i}-\frac{\partial}{\partial a_i}\right)\phi_{t-1}\Bigl(1-\su_{i=1}^n&\bigl(x_i'(a_i+1/2)+x_i''(b_i+1/2)\bigr)\Bigr).
\tag 2.17
\endalign
$$
Substituting (2.17) into (2.16) and integrating by parts we obtain
$$
\multline
\widetilde{R}_n^{(zz')}(x_1',\ldots,x_n';x_1'',\ldots,x_n'')=t^n\Gamma(t)\pro_{i=1}^n\frac{1}{x_i'-x_i''}\\
\times \inte_{\Sb a_i,b_i\\ i=1,\ldots,n\endSb}\pro_{i=1}^n\left(\frac{\partial}{\partial a_i}-\frac{\partial}{\partial b_i}\right)\frac{\phi_{-z}(a_i)\phi_{z'}(a_i+1)\phi_{-z'}(b_i)\phi_{z}(b_i+1)}{a_i+b_i+1}\\
\times \phi_{t-1}\Bigl(1-\su_{i=1}^n\bigl(x_i'(a_i+1/2)+x_i''(b_i+1/2)\bigr)\Bigr)
\pro_{i=1}^n da_idb_i.
\endmultline
$$
Simple calculations show that
$$
\multline
\left(\frac{\partial}{\partial a_i}-\frac{\partial}{\partial b_i}\right)\frac{\phi_{-z}(a_i)\phi_{z'}(a_i+1)\phi_{-z'}(b_i)\phi_{z}(b_i+1)}{a_i+b_i+1}\\
=\phi_{-z}(a_i)\phi_{z'}(a_i+1)\phi_{-z'}(b_i)\phi_{z}(b_i+1)\left(\frac{z'}{(a_i+1)b_i}-\frac{z}{a_i(b_i+1)}\right),
\endmultline
$$
and we arrive at the formula
$$
\align
\widetilde{R}_n^{(zz')}(x_1',\ldots,x_n';x_1'',\ldots,&x_n'')=t^n\Gamma(t)\pro_{i=1}^n\frac{1}{x_i'-x_i''}\\
\times \inte_{\Sb a_i,b_i\\ i=1,\ldots,n\endSb}\pro_{i=1}^n\phi_{-z}(a_i)&\phi_{z'}(a_i+1)\phi_{-z'}(b_i)\phi_{z}(b_i+1)\left(\frac{z'}{(a_i+1)b_i}-\frac{z}{a_i(b_i+1)}\right)\\
\times &\phi_{t-1}\Bigl(1-\su_{i=1}^n\bigl(x_i'(a_i+1/2)+x_i''(b_i+1/2)\bigr)\Bigr)
\pro_{i=1}^n da_idb_i.
\tag2.18
\endalign
$$
Now we have to remove all the parentheses
$$
\left(\frac{z'}{(a_i+1)b_i}-\frac{z}{a_i(b_i+1)}\right)
\tag 2.19
$$
in (2.18). We claim that the summand obtained by taking one of two terms in these parentheses for all $i=1,\ldots,n$ coincides, up to the factor 
$$
\Gamma(t)\prod\limits_{i=1}^n\phi_{z-1}(x'_i)\phi_{z'-1}(x''_i)\cdot\phi_{c-1}\Bigl(1-\frac{ |x'|+|x''|}{2}\Bigr),
$$
with the summand of (2.9) where we put $\epsilon_i=0$ if we take the first term in (2.19), and $\epsilon_i=1$ if we take the second term in (2.19). 

From this fact and (2.15) we immediately obtain the main statement (2.14) and, consequently, (2.13). 

We shall check our claim only in the case when all $\epsilon_i=1$, i.e., we take only second summands in all the parentheses (2.19). For other values of $\{\epsilon_i\}$ the considerations are quite similar.

Thus, we have to deal with the following expression
$$
\multline
t^n\Gamma(t)\pro_{i=1}^n\frac{1}{x_i'-x_i''}
 \inte_{\Sb a_i,b_i\\ i=1,\ldots,n\endSb}\pro_{i=1}^n\phi_{-z}(a_i)\phi_{z'}(a_i+1)\phi_{-z'}(b_i)\phi_{z}(b_i+1)\frac{(-z)}{a_i(b_i+1)}\\
\times \phi_{t-1}\Bigl(1-\su_{i=1}^n\bigl(x_i'(a_i+1/2)+x_i''(b_i+1/2)\bigr)\Bigr)
\pro_{i=1}^n da_idb_i.
\endmultline
$$
The  change of variables
$$
\cases
u_i=\left(1-\frac{|x'|+|x''|}{2}\right)\frac{a_i}{x_i'}\\
v_i=\left(1-\frac{|x'|+|x''|}{2}\right)\frac{b_i}{x_i''}
\endcases
$$
leads (after some simplifications) to the following expression
$$
\multline
\frac{\Gamma(t)}{(\Gamma(z)\Gamma(z'))^n}\pro_{i=1}^n\frac{(x_i')^{z'}(x_i'')^{z-1}}{x_i'-x_i''}\cdot \phi_{c-1}\Bigl(1-\frac{ |x'|+|x''|}{2}\Bigr)\\
\times \Gamma(c)\inte_{\Sb u_i,v_i\\ i=1,\ldots,n\endSb}\pro_{i=1}^n\phi_{-z-1}(u_i)\phi_{-z'}(v_i)(1-u_iy_i')^{z'}(1-v_iy_i'')^{z-1}\\
\times
 \phi_{t-1}\bigl(1-\su_{i=1}^n(u_i+v_i)\bigr)
\pro_{i=1}^n du_idv_i
\endmultline
$$
where $c=t-n(z+z'-1)$, and $y_i',y_i''$ are given by (2.14').

The last integral is the Euler-Laplace type integral representation (2.8) for $F_B^{[2n]}(a^\epsilon,b^\epsilon;c|y)$ where 
$$
\aligned
y&=(y_1',\ldots,y_n';y_1'',\ldots,y_n'')\\
a^\epsilon&=(-z',\ldots,-z';1-z,\ldots,1-z)\in \Bbb C^{2n},\\
b^\epsilon&=(-z,\ldots,-z;1-z',\ldots,1-z')\in \Bbb C^{2n},
\endaligned
$$
and the sets of parameters coincide with (2.9') for $\epsilon_1=\ldots=\epsilon_n=1$.

There is only one thing left --- to check the following equality
$$
\frac{1}{(\Gamma(z)\Gamma(z'))^n}\pro_{i=1}^n\frac{(x_i')^{z'}(x_i'')^{z-1}}{x_i'-x_i''}=
(-1)^n\pro_{i=1}^n\frac{y_i''}{y_i'-y_i''}\prod\limits_{i=1}^n\phi_{z-1}(x'_i)\phi_{z'-1}(x''_i).
$$
This can be easily done using (2.14'). \qed
\enddemo
As an application of all this heavy techniques let us consider the case $n=1$. Then Theorem 2.4.1 and Example 2.3.5 (together with (2.0)) give us the following statement.
\proclaim{Corollary 2.4.2} 
Let  $$
\aligned
a=(1-z',-z),&\qquad b=(1-z,-z');
\\
a'=(2-z',-z),&\qquad b'=(2-z,-z');
\\
a''=(1-z',1-z),&\qquad b''=(1-z,1-z');
\endaligned
$$
$$
c=-t-z-z'+1=(1-z)(1-z');
$$
$$
y=(1-1/x,1-1/x).
$$
For $x>0$ we have  
$$
\multline
\rho_1^{(zz')}(x)=\frac{\Gamma(t)}{\Gamma(z)\Gamma(z')}x^{z+z'-2}\phi_{c-1}(1-x)
 \Biggl(F_B^{[2]}(a,b;c|y)\\+y'\Bigl(F_B^{[2]}(a',b';c+1|y)-\frac{zz'}{(1-z)(1-z')}F_B^{[2]} (a'',b'';c+1|y)\Bigr)\Biggr)
\endmultline
$$
The same formula holds for $x<0$ if we replace $x$ by $|x|$ and $(z,z')$ by $(-z,-z')$.
\endproclaim
 Comparing this result with Theorem 5.12 from [O] we see that our techniques allowed to reduce the order of the multivariate hypergeometric functions involved in the expressions for the first correlation function from $3$ to $2$. Moreover, our formula allows to compute the asymptotics of the first correlation function at the origin. The result looks as follows. 
\proclaim{Corollary 2.4.3} For $x>0$
$$
\rho_1^{(zz')}(x)=\frac{A_{zz'}(x)}{x}+x^{z-z'}B_{zz'}(x)+x^{z'-z}C_{zz'}(x),
$$
where $A_{zz'}(x),B_{zz'}(x),C_{zz'}(x)$ are some functions
analytic in a neighbourhood of the origin,
$$
A_{zz'}(0)=\frac{(z-z')\sin\pi z\cdot \sin\pi z'}{\pi\sin\pi(z-z')}>0.
\tag 2.20
$$
The same formula holds for $x<0$ if we replace $x$ by $|x|$ and $(z,z')$ by $(-z,-z')$.
\endproclaim
Note that the constant $A_{zz'}(0)$ is invariant with respect to the change of $(z,z')$ by $(-z,-z')$. That means, that the density of the first controlling measure $$\sigma^{(zz')}_1(x)=|x|\rho_1^{(zz')}(x)$$ is continuous at the origin.
\demo{Sketch of the proof} If $x$ is a small positive number, then $y$ in Corollary 2.4.2 is large negative. There is a formula for the analytic continuation of $F_B^{[2]}=F_3$ to the domain of large negative values of variables which represents this function as a finite sum of expressions of the form 'minus variable raised to some (complex) power times a function analytic at the infinity', see [E, 5.11(10)]. Using this formula for all three summands in the formula of Corollary 2.4.2, we get the result. \qed
\enddemo
In Chapter 4 we shall consider much more general situation (asymptotics of the $n$th correlation function), and the $n$-dimensional variant of the proof above will be carried out in all details.

\subhead 2.5. Simplicity of the processes \endsubhead
In this section we shall prove the following statement (cf. [O, Proposition 4.2]).
\proclaim{Theorem 2.5.1} All processes $\Cal P_{zz'}$ are simple. Equivalently, for any $n\geq 2$ the set
$$
\{(x_1,\ldots,x_n)\in I^n|\ \pro_{i<j}(x_i-x_j)=0\}
$$
is a null set with respect to $\rz$.
\endproclaim
\example{Remark 2.5.2} Here is one more reformulation of the theorem (see [O, Proposition 3.5]). For almost all, with respect to $P_{zz'}$, points $(\alpha|\beta)\in \Omega$ there are no repetitions of the type $\alpha_i=\alpha_{i+1}\neq 0$ or $\beta_i=\beta_{i+1}\neq 0$.
\endexample
\example{Remark 2.5.3} Theorem 2.5.1 and Proposition 2.2.5 imply that the restriction $\prod_{i<j}(x_i-x_j)\neq 0$ in Theorems 2.2.1 and 2.4.1 can be removed (cf. Remark 2.2.2).
\endexample
\demo {Proof of Theorem 2.5.1} We shall use Proposition 3.5 from [O]. Namely, we shall prove that on the diagonal $\Delta=\{(x_1,x_2)\in I^2|\ x_1=x_2\}$ we have
$$
\sigma_2^{(zz')}|_{\Delta}=|x_1|\sigma^{(zz')}_1(dx_1)\delta(x_2-x_1).
\tag 2.22
$$
Proposition 3.5 and Proposition 4.2 of [O] state that this condition is equivalent to the statement of the theorem.

The fundamental relation (2.0) shows that it is sufficient to check (2.21) only on 
$$
\Delta_+=\{(x_1,x_2)\in I^2|\ x_1=x_2>0\}.
$$ 
 As was mentioned in Remark 1.2.2 there are three structures with two blocks. Namely,
we shall denote by $T_1$ the structure with two hook blocks, by $T_2$ --- the structure with one hook and one linear horizontal blocks, and by $T_3$ --- the structure with one hook and one linear vertical blocks. Then, see Proposition 1.3.2,
$$
\sigma^{(zz')}_2=\sigma^{(zz')}_2(T_1)+\sigma^{(zz')}_2(T_2)+\sigma^{(zz')}_2(T_3).
\tag 2.21
$$
We start with $T_2$. Let us compute $\sigma^{(zz')}_2(T_2)$ using Proposition 1.5.2. We shall omit the subscript $i$ because our structure contains only one fragment, and $i\equiv 1$. We get
$$
\sigma^{(zz')}_2(T_2)=t\Cal C^a_{a'a''}\Cal D^{a'b}_{a'}\left[a'I_1(a',a'')
\right]
$$
where
$$
x_1=a,\quad x_2=b.
$$
Applying Propositions 1.4.1 and 1.4.3 we obtain
$$
\multline
\sigma^{(zz')}_2(T_2)=t\Cal C^a_{a'a''}\left[a'I_1(a',a'')\cdot\delta(b-a')\right]\\
=ta^2\inte_yI_1(a,y)\frac{dy}{a+y}\cdot\delta(b-a)+\frac{t|a|b}{|a|+b}I_1(b,-a).
\endmultline
$$
The relation (1.7) implies that the second summand does not contribute to the positive quadrant.
Moreover, in the first summand we recognize  the first controlling measure multiplied by $a\cdot\delta(b-a)=x_1\delta(x_2-x_1)$, cf. (2.1). Comparing this to (2.21) we see, that we have to prove that $\Delta_+$ ia a null set for $\sigma^{(zz')}_2(T_1)+\sigma^{(zz')}_2(T_3)$.

Arguing as above, using Proposition 1.5.2, one can show that $\sigma^{(zz')}_2(T_3)$ gives no contribution to the positive quadrant. Furthermore, from Proposition 1.5.2 it is easy to derive that the contribution of $\sigma^{(zz')}_2(T_1)$ to the positive quadrant  equals (cf. (2.1))
$$
t^2x_1x_2\inte_{s_1,s_2}\frac{I_2(x_1,s_1;x_2,s_2)}{(x_1+s_1)(x_2+s_2)}ds_1ds_2.
$$
In 2.2 (Theorem 2.2.1) we proved that this is the same as the right-hand side of (2.2) taken for $n=2$, multiplied by  $x_1x_2$. 
 Then, by Proposition 2.2.5, it is an analytic function in the domain
$$
\{(x_1,x_2)\in I^2|\ x_1,x_2>0,\ x_1+x_2\neq1\}.
$$
(Easy to see from (2.1), it is identically zero if $x_1+x_2>1$.) That is why the only possibility for $\sigma^{(zz')}_2(T_1)$ to have $\Delta_+$ as a non-null set is to have an atom at the point 
$$
(x_1,x_2)=(1/2,1/2).
$$
If it were the case, then $\sigma_2^{(zz')}$ would have an atom at the same point, and
$$
\sigma_1^{(zz')}(dx_1)=\int_{x_2}\sigma_2^{(zz')}(dx_1,dx_2)
$$
would have an atom at the point $x_1=1/2$. This conclusion contradicts with Corollary 2.2.7.\qed
\enddemo 
From Theorems 2.2.1, 2.5.1, and Corollary 2.2.5 we immediately obtain (cf. Corollary 2.2.7)
\proclaim{Corollary 2.5.4} For any $n=1,2,\ldots$ the correlation functions $\rz$ are analytic in the domain   
$$
\{(x_1,\ldots,x_n)|\ x_1,\ldots,x_n>0,\ \sum_{i=1}^n x_i <1\}.
$$
\endproclaim
Note that all the correlation functions are identically zero as soon as $\sum|x_i|>1$.
\example{Remark 2.5.5} We suppose that by explicit formulas of Chapter 1 a more general statement can be proved. Namely, that for all $n=1,2,\ldots$ the $n$th correlation function is analytic in the domain
$$
\{(x_1,\ldots,x_n)\in I^n|\ \sum_{i=1}^n|x_i|<1\}.
$$
\endexample 
\head 3. Lifted Processes \endhead
\subhead 3.1. Lifting. Lifted correlation functions \endsubhead
 Two previous chapters were devoted to the point processes
 originated from the measures $P_{zz'}$ on the Thoma simplex 
 $\Omega$ as was explained in Introduction (see also [O]). In fact,
 G.~I.~Olshanski proved, see [O, Theorem 6.1], that all these measures are concentrated on the face $\Omega_0$ of $\Omega$. Let us define a covering $\widetilde{\Omega_0}$ of $\Omega_0$ as follows 
$$
\widetilde{\Omega_0}=\Omega_0\times \Bbb R_+
$$
(we understand this product as the product of Borel spaces).

There is a natural one-to-one correspondence between $\widetilde{\Omega_0}$ and the set
$$
\widetilde{\Omega_0'}= \{\tilde{\alpha_1}\geq\tilde{\alpha_2}\geq 0,\
\tilde{\beta_1}\geq\tilde{\beta_2}\geq 0; \ \sum\limits_i(\tilde{\alpha_i}+\tilde{\beta_i})=s\in \Bbb R_+\}.
$$
Namely, the map 
$$
f:\widetilde{\Omega_0}\to \widetilde{\Omega_0'};\quad
f:((\alpha|\beta),s)\mapsto (s\alpha|s\beta)
$$
is bijective. In the following we identify $\widetilde{\Omega_0}$ and $\widetilde{\Omega_0'}$.

Starting from any probability measure $\mu$ on $\Omega_0$ let us introduce the probability measure $\tilde{\mu}$ on $\widetilde{\Omega_0}$ as follows
$$
\tilde{\mu}=\mu\otimes \frac{s^{t-1}}{\Gamma(t)}e^{-s}ds.
$$

After a while we shall set $\mu=P_{zz'}$, but the general construction is of certain independent interest (see, for example, Proposition 3.1.2).

Following the general scheme described in [O, \S 4], we associate with $\mu$ a stochastic point process $\Cal P_{\mu}$ on $I=[-1,1]\setminus\{0\}$. Let us denote its correlation functions by $\rho_n(x_1,\ldots,x_n)$.
 
Now we shall construct a new stochastic point process. Consider the set $\widetilde{I}=\Bbb R\setminus \{0\}$ as the phase space; sets that do not intersect a sufficiently small interval $(-\varepsilon,\varepsilon)$ as the test sets; finite and countable systems of points in $\widetilde{I}$, such that their intersection with any test set is finite, as the configurations.

To any point $(\tilde{\alpha}|\tilde{\beta})\in \widetilde{\Omega_0'}=\widetilde{\Omega_0}$
we attach the configuration $$(\tilde{\alpha_1},\tilde{\alpha_2},\ldots,
-\tilde{\beta_1},-\tilde{\beta_2},\ldots).$$ Then our measure $\tilde\mu$ makes this configuration random. We are interested in the correlation functions of the point process $\widetilde{\Cal P_\mu}$ thus obtained. We shall call $\widetilde{\Cal P_\mu}$ the {\it lifting} of the initial process $\Cal P_\mu$.
\proclaim{Proposition 3.1.1} All correlation measures $\tilde{\rho}_n$ of $\widetilde{\Cal P_\mu}$ take finite values on test sets. Furthermore,
$$
\tilde{\rho}_n(x_1,\ldots,x_n)=\int\limits_0^{+\infty} \frac {s^{t-1}}{\Gamma(t)}e^{-s}\rho_n(x_1s^{-1},\ldots,x_ns^{-1})\frac{ds}{s^n}.
\tag 3.1
$$
\endproclaim
Under appropriate assumptions about $\rho_n$, this integral transform can be inverted.  

\demo{Proof (due to G.~I.~Olshanski)} Take a test set $A$ of $\Cal P_\mu$ and consider
 a random variable $N_A$ on $\Omega_0$ defined as follows. The value of $N_A$ on 
 $(\alpha|\beta)\in \Omega_0$ is the number of $\alpha_i$'s and $-\beta_i$'s which 
 belong to $A$. By [O, Proposition 4.1] $N_A$ is everywhere finite, and for
$$
A=A(\varepsilon)=[-1,-\varepsilon)\cup (\epsilon,1]
$$ 
we have
$$
N_{A(\varepsilon)}\leq \frac{1}{\varepsilon}.
$$

For a test $\widetilde A$ of $\widetilde {\Cal P_{\mu}}$ we define the random variable
 $\widetilde N_{\widetilde A}$ in the same way: it counts the number of $\tilde\alpha_i$'s
  and $-\tilde\beta_i$'s in $\widetilde A$.

By definition of correlation measures, the value of the $n$th correlation measure $\tilde\rho_n$ on the test set $\widetilde{A}^n$ is the factorial moment of $\widetilde N_{\widetilde A}$:
$$
\tilde\rho_n(\widetilde{A}^n)=\Bbb E(\widetilde N_{\widetilde A}(\widetilde N_{\widetilde A}-1)\cdots(\widetilde N_{\widetilde A}-n+1)).
$$
Thus, it suffices to prove, that the usual moments of  $\widetilde N_{\widetilde A}$ are finite. We may consider only test sets of the form
$$
\widetilde A=\widetilde A(\varepsilon)=(-\infty,-\varepsilon)\cup (\epsilon,+\infty).
$$ 
For any $(\alpha|\beta)\in \Omega_0$ and $s>0$ we have
$$
\widetilde N_{\widetilde A(\varepsilon)}(s\alpha|s\beta)=N_{A(\varepsilon s^{-1})}(\alpha|\beta)\leq \frac{s}{\varepsilon}.
$$
Hence,
$$
\multline
\Bbb E(\widetilde N_{\widetilde A(\varepsilon)}^k)=\inte_{\widetilde{\Omega_0}}\widetilde N_{\widetilde A(\varepsilon)}^k
\tilde{\mu}(d\tilde\omega)=\inte_0^{\infty}\inte_{\Omega_0} N_{ A(\varepsilon s^{-1})}^k\mu(d\omega)\frac{s^{t-1}}{\Gamma(t)}e^{-s}ds\\
\leq\frac{1}{\varepsilon^k}\inte_0^{\infty}\frac{s^{k+t-1}}{\Gamma(t)}e^{-s}ds=\frac{(t)_k}{\varepsilon^k}<\infty. 
\endmultline
$$
The formula (3.1) is obvious from the definition of the correlation functions. \qed
\enddemo

In the following sections we shall see that the lifting substantially simplifies the formulas
 for the correlation functions of our processes $\Cal P_{zz'}$. Moreover, it also works for 
 the Poisson-Dirichlet processes, see [Ki] for definitions.
 
 The following claim also follows from section 9.4 in [Ki].
 
\proclaim{Proposition 3.1.2} The lifting of the Poisson-Dirichlet process $\Cal {PD}(t)$ is the Poisson process on $(0,+\infty)$ with density $te^{-x}/x$.
\endproclaim
\demo{Proof}
Simple calculation. The correlation functions $\Cal {PD}(t)$ have the form, see [O, Corollary 7.4],
$$
\rho_n(x_1,\ldots,x_n)=\frac{t^n(1-x_1-\ldots-x_n)^{t-1}_+}{x_1\cdots x_n}.
$$
Then by (3.1) the lifted correlation functions are
$$
\multline
\tilde{\rho}_n(x_1,\ldots,x_n)=\int\limits_0^{+\infty} \frac {s^{t-1}}{\Gamma(t)}e^{-s}\frac{t^n(1-x_1s^{-1}-\ldots-x_ns^{-1})^{t-1}_+}{x_1\cdots x_n}
ds\\=\frac{t^n e^{-x_1-\ldots-x_n}}{x_1\cdots x_n}.\qed
\endmultline
$$
\enddemo

\subhead 3.2. General structure of  $\tilde{\rho}_n^{(zz')}$
\endsubhead
From now on we apply the lifting to $\mu=P_{zz'}$.
Our main goal in this section is to reformulate Theorem 1.7.3 in terms of lifted correlation functions. To do this we need to extend the notion of lifting to the distributions with compact supports. In what follows $\tau>0$ is a positive number.
\proclaim{Proposition 3.2.1} Let $f(\xi_1,\ldots,\xi_k)$ be a distribution in $k$ variables with  compact support. Then there exists a distribution
$$
(\Cal L^{\tau} f)(\xi_1,\ldots,\xi_k)=\inte_0^{\infty} f\left(\frac{\xi_1}{s},\ldots,\frac{\xi_k}{s}\right)\frac{s^{\tau-k-1}}{\Gamma(\tau)}e^{-s}{ds}
\tag 3.2
$$
not necessarily with a compact support, such that its values on a compactly supported test function $\psi(\xi_1,\ldots,\xi_k)$ is given by the formula
$$
\langle \tilde f,\psi\rangle=\inte_0^{\infty}\frac{s^{\tau-1}}{\Gamma(\tau)}e^{-s}\left[\inte_{\xi_1,\ldots,\xi_k}
f(\xi_1,\ldots,\xi_k)\psi(s\xi_1,\ldots,s\xi_k)d\xi_1,\cdots,d\xi_k\right] ds.
\tag 3.3
$$
\endproclaim
\example{Remark 3.2.2} Formulas (3.1) and (3.2) represent the same operation if $\tau=t$.
\endexample
\demo{Proof of Proposition 3.2.1} The only thing we need to check is the convergence of the integral over $s$ in (3.3). But this becomes obvious if we remember that any distribution with  compact support is the result of applying some differential operator $\partial ^{m_1}_{\xi_1} \cdots\partial ^{m_k}_{\xi_k}$ to a continuous function with  compact support.\qed
\enddemo
Note that (3.3) is correctly defined if $\psi$ is any bounded test function (not necessarily with compact support).

Now we want to extend the notion of pseudoconvolution, see 1.4, to distributions which do not necessarily have compact supports. Namely, we shall say that a distribution $h(\xi_1,\ldots,\xi_k)$ is the pseudoconvolution of distributions $f(\xi_1,\ldots,\xi_k)$ and
$g(\xi_1,\ldots,\xi_k)$ (and write $h=f\odot g$) if for any compactly supported test function $\psi(\xi_1,\ldots,\xi_k)$ the following relation holds
$$
\langle h,\psi\rangle=\int f(\xi_1,\ldots,\xi_m)g(\eta_1,\ldots,\eta_m)\psi(\xi_1\eta_1,\ldots,\xi_m\eta_m)\pro_id\xi_i d\eta_i
$$
\proclaim{Proposition 3.2.3} For any distributions
$f(\xi_1,\ldots,\xi_k)$ and $g(\xi_1,\ldots,\xi_k)$ with \linebreak compact supports
$$
\Cal L^{\tau}({f\odot g})=(\Cal L^{\tau} f)\odot g.
\tag 3.4
$$
\endproclaim
\demo{Proof} Obvious.
\enddemo
\proclaim{Proposition 3.2.4}
The lifted correlation function have the form
$$
\tilde\rho^{(zz')}_n(x_1\ldots,x_n)=\sum\limits_{d\geq n/2}^n \sum \limits_{\varphi \in \Phi_{n,d}}(\varphi \widetilde H_d)(x_1,\ldots,x_n)
$$
where
$$
\widetilde H_d(r_1,s_1;\ldots;r_d,s_d)=\frac{t^d(t)_d}{d!\prod\limits_{i=1}^d(r_i+|s_i|)} (\Cal L^{t+d}I_d)(r_1,-s_1;\ldots;r_d,-s_d).
$$
\endproclaim
\demo{Proof}
Combining (3.1) and Theorem 1.7.3 we see that it is sufficient to prove the following relation
$$
\Cal L^{t}(\varphi H_d)=\frac{\Gamma(t+d)}{\Gamma(t)}\varphi \widetilde{H_d}
\tag 3.5
$$
for all $\varphi\in \Phi_{n,d}$. Let us check this for one particular $\varphi$, for all others the proof is quite similar. Let us take 
$$
\varphi:\{1,\ldots,n\}\to \{1,1';\ldots;d,d'\}
$$
such that
$$
\align
\varphi(i)=i',&\quad i\leq d  \\
\varphi(i)=i-d,&\quad i>d.
\endalign
$$
Then, see 1.7,
$$
\multline
(\varphi H_d)(x_1,\ldots,x_n)=\frac{t^d}{\prod_{i=1}^{n-d}(x_{i+d}+|x_i|)}\inte_{r_{n-d+1},\ldots,r_d}
\frac{1}{\prod_{i=n-d+1}^d(r_i+|x_i|)}\\
\times
I_d(x_{d+1},-x_1;\ldots;x_{n},-x_{n-d};r_{n-d+1},-x_{n-d+1};\ldots;r_d,-x_d)\pro_{i=n-d+1}^ddr_i.
\endmultline
$$
Hence, 
$$
\multline
\Cal L^{t}(\varphi H_d)(x_1,\ldots,x_n)=\inte_0^{\infty}\frac{s^{t-n-1}}{\Gamma(t)}e^{-s}\\
\times \frac{t^d}{\prod_{i=1}^{n-d}(x_{i+d}/s+|x_i/s|)}\inte_{r_{n-d+1},\ldots,r_d}
\frac{1}{\prod_{i=n-d+1}^d(r_i+|x_i/s|)}
\\
\times I_d\left(\frac{x_{d+1}}{s},\frac{-x_1}{s};\ldots;\frac{x_{n}}{s},\frac{-x_{n-d}}{s};r_{n-d+1},\frac{-x_{n-d+1}}{s};\ldots;r_d,\frac{-x_d}{s}\right)\pro_{i} dr_i\cdot ds.
\endmultline
$$
By changing the variables 
$$
u_i=sr_i,\quad i=n-d+1,\ldots,d,
$$
we arrive at the formula
$$
\multline
\frac{(-1)^dt^d}{\prod_{i=1}^{n-d}(x_{i+d}+|x_i|)}\inte_{u_{n-d+1},\ldots,u_d}
\frac{1}{\prod_{i=n-d+1}^d(u_i+|x_i|)}\inte_0^{\infty}\frac{s^{t-d-1}}{\Gamma(t)}e^{-s}\\
\times I_d\left(\frac{x_{d+1}}{s},\frac{-x_1}{s};\ldots;\frac{x_{n}}{s},\frac{-x_{n-d}}{s};\frac{u_{n-d+1}}{s},\frac{-x_{n-d+1}}{s};\ldots;\frac{u_d}{s},\frac{-x_d}{s}\right)ds\cdot \pro_{i} du_i
\endmultline
$$
that coincides with the right-hand side of  (3.5). \qed
\enddemo
It remains to compute $(\Cal L^{t+d}I_d)$.
\proclaim{Lemma 3.2.5} Let $\alpha_1+\ldots+\alpha_m=0$. Then
$$
\Cal L^{\alpha_0+m+1}\Bigr(\phi_{\alpha_1}(u_1)\cdots\phi_{\alpha_m}(u_m)\phi_{\alpha_0}\bigl(1-\sum_{i=1}^mu_i\bigl)\Bigl)=\frac{\phi_{\alpha_1}(u_1)\cdots\phi_{\alpha_m}(u_m)e^{-\sum_{i=1}^mu_i}}{{\Gamma(\alpha_0+m+1)}}.
 $$
\endproclaim
\demo{Proof} Direct computation.
$$
\multline
\inte_0^{\infty}\frac{s^{\alpha_0}}{\Gamma(\alpha_0+m+1)}e^{-s}\phi_{\alpha_1}(u_1/s)\cdots\phi_{\alpha_m}(u_m/s)\phi_{\alpha_0}\bigl(1-\sum_{i=1}^mu_i/s\bigl)ds\\
=\inte_0^{\infty}\frac{s^{-\sum_{i=1}^m\alpha_i}}{\Gamma(\alpha_0+m+1)}e^{-s}\phi_{\alpha_0}\bigl(s-\sum_{i=1}^mu_i\bigl)ds\cdot \phi_{\alpha_1}(u_1)\cdots\phi_{\alpha_m}(u_m)\\
=\frac{1}{\Gamma(\alpha_0+m+1)}\inte_{0}^{\infty}\phi_{\alpha_0}\bigl(s)e^{s-\sum_{i=1}^mu_i}ds\cdot \phi_{\alpha_1}(u_1)\cdots\phi_{\alpha_m}(u_m)\\
=\frac{\phi_{\alpha_1}(u_1)\cdots\phi_{\alpha_m}(u_m)e^{-\sum_{i=1}^mu_i}}{\Gamma(\alpha_0+m+1)}.\qed
\endmultline
$$
\enddemo
Let us use the explicit formula (1.22) for  $I_d$. The previous lemma shows that
$$
\Cal L^{t+d}F_1(a_1',a_1'';\ldots;a_d',a_d'')=\frac{\Gamma(t)}{\Gamma(t+d)}\pro_{i=1}^d \phi_{z'}(a_i')\phi_{-z'}(a_i'')e^{-a_i'-a_i''}. 
$$
Thus, using (3.4), we obtain
$$
(\Cal L^{t+d}I_d)(a_1',a_1'';\ldots;a_d',a_d'')=\frac{\Gamma(t)}{\Gamma(t+d)}\det N(a_i',a_i'')
$$
where
$$
N(a,b)=\left(\phi_{z'}(a)\phi_{-z'}(b)e^{-a-b}\right)\odot \left(\phi_{z,-z-1}(a)\phi_{-z,z-1}(b)\right)\odot\left(\delta(a-b)\chi_{[0,1]}(a)\right).
\tag 3.6
$$
This gives us, together with Proposition 3.2.4, the main statement of this section.
\proclaim{Theorem 3.2.6}
$$
\tilde{\rho}_n^{(zz')}(x_1,\ldots,x_n)=
\sum\limits_{d\geq n/2}^{n} \sum\limits_{\varphi\in \Phi_{n,d}}(\varphi \widetilde{H_d})(x_1,\ldots,x_n)
$$
where
$$
\widetilde{H_d}(r_1,s_1;\ldots;r_d,s_d)=\frac{t^d}{d!\prod\limits_i (r_i+|s_i|)}\det N(r_i,-s_j);
$$
and $N(a,b)$ is given by (3.6).
\endproclaim
\example{Remark 3.2.7}
Clearly,
$$
\operatorname{supp}H_d \subset \{(r_1,s_1;\ldots;r_d,s_d)|\ r_i\geq 0,s_i\leq 0,i=1,\ldots,d\}.
$$
\endexample
\example{Remark 3.2.8} Our main achievement is that the distributions $I_d$ which are rather
complicated functions in $2d$ variables, were `reduced' to one function $N$ in two variables. 
\endexample  
\subhead 3.3. Lifted correlation functions in positive (negative) hyperoctants \endsubhead
As was mentioned before (the beginning of Chapter 2) we may consider only positive hyperoctants. By analogy with Proposition 2.1.1 we immediately see, using Theorem 3.2.6, that in the positive hyperoctant
$$
\tilde \rz(x_1,\ldots,x_n)=\det M(x_i,x_j)
$$
where
$$
M(x,y)=\inte_s\frac{tN(x,s)}{s+y}ds
$$
is the Stieltjes transform of the kernel $tN(r,s)$ with respect to the second argument. 
In fact, we have already computed the kernel $M$. Namely, we shall derive from Theorem 2.2.1 the following statement.
\proclaim{Theorem 3.3.1}
Let $x_1,\ldots,x_n>0$. Then
$$
\tilde{\rho}_n^{(zz')}(x_1,\ldots,x_n)=\det M(x_i,x_j)
\tag3.7
$$
where 
$$
\aligned
M(x,y)=t\iint\limits_{t_1,t_2}&\phi_{-z}(t_1)\phi_{-z'}(t_2)\phi_{z'}(t_1+1)\phi_{z}(t_2+1)\frac{e^{-x(t_1+1/2)-y(t_2+1/2)}}{t_1+t_2+1}dt_1dt_2.
\endaligned
$$
\endproclaim
\demo{Proof}
We shall use the fact, see [KOV], that the set
$$
\{(x_1,\ldots,x_n)\in I^n|\ \sum_{i=1}^n|x_i|=1\}
$$
is a null set for the $n$th correlation function $\rz$.  That is why we shall neglect this set  while making the lifting (3.1).

Let us apply $\Cal L^{t}$ to both sides of (2.2). Using the relation
$$
\Cal L^t\left[\phi_{t-n-1}\bigl(1-\su_{i=1}^nx_i(a_i+b_i+1)\bigr)\right]=\frac{1}{\Gamma(t)}{e^{-\sum_{i=1}^n x_i(a_i+b_i+1)}}
$$
we arrive at our assertion. \qed
\enddemo
\example{Remark 3.3.2} An immediate corollary of the last theorem is that the lifted correlation functions are analytic in the positive hyperoctants. We suppose that using Theorem 3.2.6 one can prove that lifted correlation functions are analytic whenever all variables are nonzero (cf Remark 2.5.5). Note that the initial correlation functions were not even continuous for some values of $z$ and $z'$ on the borders $\sum_{i=1}^n|x_i|=1$ of their supports. (For the first correlation function it can be seen from the explicit formula in Corollary 2.4.2).
\endexample
\example{Remark 3.3.3} Theorem 3.3.1 means that if we restrict ourselves to the behaviour of $\{\tilde\alpha_i\}$ only (i.e. we consider the point process in the restricted phase space $\Bbb R_+\subset\widetilde{I}$) then all correlation functions are given by the determinants of the type (3.7).
Such processes were considered by several authors, see [Me], [DVJ]. Furthermore, determinantal formulas for correlation functions appear in some models of mathematical physics,  see [KBI]. For a more detailed discussion of these processes see [BO].
\endexample

  It turns out that the main result of Section 2.4 (Theorem 2.4.1) also has a nice reformulation in terms of the lifted processes. It expresses the kernel $M(x,y)$ via the Whittaker function $W_{\kappa,\mu}(x)$, see [E, chapter 6] for the definition. 
\proclaim{Theorem 3.3.4}
Let $x_1,\ldots,x_n>0$. Then
$$
\tilde{\rho}_n^{(zz')}(x_1,\ldots,x_n)=\det K(x_i,x_j)
$$
where
$$
K(x,y)=\frac{1}{\Gamma(z)\Gamma(z')}\frac{\varphi_1(x)\varphi_2(y)-\varphi_1(y)\varphi_2(x)}{x-y};
\tag 3.8
$$
$$
\varphi_1(x)=x^{-\frac{1}{2}}W_{\frac{z+z'+1}{2},\frac{z-z'}{2}}(x),\quad 
\varphi_2(x)=x^{-\frac{1}{2}}W_{\frac{z+z'-1}{2},\frac{z-z'}{2}}(x),
$$
and $W_{\kappa,\mu}(x)=W_{\kappa,-\mu}(x)$ is the Whittaker function.
\endproclaim

Clearly, the kernel $K(x,y)$ is real symmetric:
$$
K(x,y)=K(y,x).
$$
We call $K(x,y)$ the {\it Whittaker kernel}.
\demo{Proof}
We shall prove the formula
$$
M(x,y)=(x/y)^{\frac{z-z'}{2}}K(x,y),
\tag 3.9
$$
which by Theorem 3.3.1 proves the assertion.

The simplest way to prove (3.9) is to apply $\Cal L^{t+n}$ to both sides of (2.18). However, we prefer to give an independent proof here.

We start from the integral representation of $M(x,y)$ in Theorem 3.3.1. Note that
$$
(x-y)e^{-x(t_1+1/2)-y(t_2+1/2)}=\left(\frac{\partial}{\partial t_2}-\frac{\partial}{\partial t_1}\right)e^{-x(t_1+1/2)-y(t_2+1/2)}.
$$
Integration by parts gives
$$
\multline
(x-y)M(x,y)=
t\iint\limits_{t_1,t_2}e^{-x(t_1+1/2)-y(t_2+1/2)}\\ \times
\left(\frac{\partial}{\partial t_1}-\frac{\partial}{\partial t_2}\right)\left[\frac{\phi_{-z}(t_1)\phi_{-z'}(t_2)\phi_{z'}(t_1+1)\phi_{z}(t_2+1)}{t_1+t_2+1}\right]dt_1dt_2.
\endmultline
$$
Simple calculation shows that
$$
\multline
\left(\frac{\partial}{\partial t_1}-\frac{\partial}{\partial t_2}\right)\frac{\phi_{-z}(t_1)\phi_{-z'}(t_2)\phi_{z'}(t_1+1)\phi_{z}(t_2+1)}{t_1+t_2+1}\\
=\phi_{-z}(t_1)\phi_{-z'}(t_2)\phi_{z'}(t_1+1)\phi_{z}(t_2+1)\left(\frac{z'}{(t_1+1)t_2}-\frac{z}{t_1(t_2+1)}\right).
\endmultline
$$
Then we get
$$
\multline
(x-y)M(x,y)\\=\frac{e^{\frac{-x-y}{2}}}{\Gamma(z)\Gamma(z')}\Biggl(\inte_{t_1}\phi_{-z-1}(t_1)(t_1+1)^{z'}e^{-xt_1}dt_1\inte_{t_2}\phi_{-z'}(t_2)(t_2+1)^{z-1}e^{-yt_2}dt_2\\ -
\inte_{t_1}\phi_{-z}(t_1)(t_1+1)^{z'-1}e^{-xt_1}dt_1\inte_{t_2}\phi_{-z'-1}(t_2)(t_2+1)^{z}e^{-yt_2}dt_2\Biggr).
\endmultline
$$
Using the standard integral representation of the Whittaker function (or that of the confluent hypergeometric function $\Psi(a,c;x)$, see [E, 6.9.4 and 6.5.2])
$$
W_{\kappa,\mu}(x)=e^{-x/2}x^{\mu+1/2}\inte_{\tau}\phi_{\mu-\kappa-1/2}(\tau)(1+\tau)^{\mu+\kappa-1/2}e^{-\tau x}d\tau
$$
we arrive at (3.9). \qed
\enddemo
\head 4. Asymptotics at the Origin\endhead
  In this chapter we compute the asymptotics of both lifted and non-lifted correlation functions
  when all variables are of the same sign and infinitely small. The results can be considered as
  applications of the explicit formulas obtained in the previous chapters. The asymptotics is
  also of certain independent interest; it gives a possibility to consider new (stationary)
  stochastic processes obtained from the points which are `infinitely close to zero', see [BO] for details. In other words, these new processes encode the information about the behaviour of $\alpha_k, \beta_k$ as $k\to \infty$.

The idea of considering the asymptotics of the correlation
functions is due to G.~I.~Olshanski. He also computed the asymptotics in the lifted case (Section 4.1), and this served as a prompt for our further results in this chapter.

\subhead 4.1. Lifted processes \endsubhead
In this section we shall prove the following statement.
\proclaim{Theorem 4.1.1} Let $x_1,\ldots,x_n>0$. Then
$$
\tilde \rho_n^{(zz')}(x_1,\ldots,x_n)=
\frac{\det k({x_i}/{x_j})
+r(x_1,\ldots,x_n)}{x_1\cdots x_n}
\tag 4.1
$$
where
$$
k(x)=\frac{\sin\pi z\cdot \sin\pi z'}{\pi\sin\pi(z-z')}
\cdot
\frac{x^{\frac{z-z'}2}-x^{\frac{z'-z}2}}{x^{\frac12}-x^{-\frac12}};
\tag 4.2
$$
$$
r(x_1,\ldots,x_n)=\cases
O\left(\max\{x_i\}\right),\ &\text{if}\ z'=\bar z\ne z\\
O\left(\max\{x_i\ln^2 x_i\}\right),\ 
&\text{if}\ z'=z\in\Bbb R\setminus\Bbb Z\\
O\left({(\max\{x_i\})}^{1-|z-z'|}\right),\ 
&\text{if}\ m<z,z'<m+1,\ m\in\Bbb Z
\endcases
$$
as $x_1,\ldots,x_n\to 0$.
\endproclaim
\example{Remark 4.1.2} The main term of the asymptotics depends 
only on the ratios of the variables. When $z'=z$, the function $k(x)$
is defined as the limit of the expression \thetag{4.2} as $z'-z\to0$.
\endexample 
\demo{Proof}
By the determinantal formula for lifted correlation functions 
proved in Theorem 3.3.4, it suffices to compute the asymptotics 
of the Whittaker kernel $K(x,y)$. It turns out that the 
cases $z\neq z'$ and $z=z'$ require different treating. The latter 
is the limit case of the former, and the asymptotics
for $z=z'$  (Proposition 4.1.4) can be formally obtained from 
that for $z\neq z'$ (Proposition 4.1.3).

The claim of our theorem immediately follows from two following statements.
\proclaim{Proposition 4.1.3} Let $x,y>0$ and $z\neq z'$. Then
$$
\gather
K(x,y)=\frac{\sin\pi z\cdot \sin\pi z'}{\pi\sin\pi(z-z')}
\cdot\frac{(x/y)^{\frac{z-z'}{2}}-(x/y)^{\frac{z'-z}{2}}}{x-y}
+\frac{r'(x,y)}{\sqrt{xy}}\\
=\frac1{\sqrt{xy}}\left\lbrace
\frac{\sin\pi z\cdot \sin\pi z'}{\pi\sin\pi(z-z')}
\cdot\frac{(x/y)^{\frac{z-z'}{2}}-(x/y)^{\frac{z'-z}{2}}}
{(x/y)^{\frac12}-(x/y)^{-\frac12}}
+r'(x,y)
\right\rbrace
\endgather
$$
where
$$
r'(x,y)=\cases
O\left(\max\{x,y\}\right),\ &\text{if}\ z'=\bar z\ne z\\
O\left({(\max\{x,y\})}^{1-|z-z'|}\right),\ 
&\text{if}\ m<z,z'<m+1,\ m\in\Bbb Z
\endcases
$$
as $x,y\to 0$
\endproclaim
\proclaim{Proposition 4.1.4} Let $x,y>0$ and $z=z'$. Then
$$
\gather
K(x,y)=\frac{\sin^2\pi z}{\pi^2}\cdot\frac{\ln x-\ln y}{x-y}
+\frac{r''(x,y)}{\sqrt{xy}}\\
=\frac1{\sqrt{xy}}\left\lbrace
\frac{\sin^2\pi z}{\pi^2}\cdot\frac{\ln x-\ln y}
{(x/y)^{\frac12}-(x/y)^{-\frac12}}+r''(x,y)
\right\rbrace
\endgather
$$
where
$$
r''(x,y)=
O\left(\max\{x\ln^2 x,y\ln^2 y\}\right)
$$
as $x,y\to 0$
\endproclaim
In what follows we shall use the symbol $\approx$ to relate
expressions the difference of which has the same order as $r'(x,y)$
or $r''(x,y)$ (depending on the case we consider).  

\demo{Proof of Proposition 4.1.3}
We start with the following formula that follows from the definition 
of the Whittaker function and basic properties of confluent 
hypergeometric functions, see [E], 6.5(6) and 6.9(2). 
For $\mu\notin \frac{1}{2}\Bbb Z$
$$
\multline
x^{-{1}/{2}}W_{\kappa,\mu}(x)
=e^{-{x}/{2}}\Bigl(\frac{\Gamma(-2\mu)x^{\mu}}
{\Gamma(1/2-\kappa-\mu)}\Phi(1/2-\kappa+\mu,2\mu+1;x)\\
+ \frac{\Gamma(2\mu)x^{-\mu}}
{\Gamma(1/2-\kappa+\mu)}\Phi(1/2-\kappa-\mu,-2\mu+1;x)\Bigr)
\endmultline
$$
where
$$
\Phi(a,c;x)=\sum_{k=0}^\infty\frac{(a)_k}{(c)_k}\frac{x^k}{k!}
$$
is the confluent hypergeometric function. We need only the first
terms of each summand of the expression above. We get
$$
\varphi_1(x)=x^{\frac{z-z'}{2}}A_1(x)+x^{\frac{z'-z}{2}}B_1(x)
$$
$$
\varphi_2(x)=x^{\frac{z-z'}{2}}A_2(x)+x^{\frac{z'-z}{2}}B_2(x)
$$
where the functions $A_i(x),B_i(x),\  i=1,2,$ are analytic in a 
neighbourhood of zero and
$$
A_1(0)=\frac{\Gamma(z'-z)}{\Gamma(-z)},\quad B_1(0)
=\frac{\Gamma(z-z')}{\Gamma(-z')};
$$
$$
A_2(0)=\frac{\Gamma(z'-z)}{\Gamma(1-z)},\quad B_2(0)
=\frac{\Gamma(z-z')}{\Gamma(1-z')}.
$$
Then
$$
\gather
K(x,y)=\frac{1}{\Gamma(z)\Gamma(z')}
\Biggl((xy)^{\frac{z-z'}{2}}\frac{A_1(x)A_2(y)-A_1(y)A_2(x)}{x-y}\\
+(xy)^{\frac{z'-z}{2}}\frac{B_1(x)B_2(y)-B_1(y)B_2(x)}{x-y}\\
+\frac{(x/y)^{\frac{z-z'}{2}}
A_1(x)B_2(y)-(x/y)^{\frac{z'-z}{2}}A_1(y)B_2(x)}{x-y}\\
+\frac{(x/y)^{\frac{z'-z}{2}}
B_1(x)A_2(y)-(x/y)^{\frac{z-z'}{2}}B_1(y)A_2(x)}{x-y}\Biggr)
\endgather
$$
This formula consists of four summands, let us denote them 
by $S_1,S_2,S_3,S_4$ respectively.

Note that
$$
(xy)^{\frac12}x^{\pm\frac{z-z'}2}y^{\mp\frac{z-z'}2}
=x^{\frac12\pm\frac{z-z'}2}y^{\frac12\mp\frac{z-z'}2}\approx0.
$$
Indeed, if $z'=\bar z$ then 
$$
\left|x^{\frac12\pm\frac{z-z'}2}y^{\frac12\mp\frac{z-z'}2}\right|
=x^{\frac12}y^{\frac12}\le\max\{x,y\}\approx0,
$$
and if $m<z,z'<m+1$ for a certain $m\in\Bbb Z$ then
$|\frac{z-z'}2|<\frac12$ and
$$
0<x^{\frac12\pm\frac{z-z'}2}y^{\frac12\mp\frac{z-z'}2}
\le x^{\frac12-|\frac{z-z'}2|}y^{\frac12-|\frac{z-z'}2|}
\le\left(\max\{x,y\}\right)^{1-|z-z'|}\approx 0.
$$

Further, the functions 
$$
\frac{A_1(x)A_2(y)-A_1(y)A_2(x)}{x-y}\,,\qquad
\frac{B_1(x)B_2(y)-B_1(y)B_2(x)}{x-y}
$$
are analytic and consequently bounded in a neighbourhood of the
origin. Hence, $(xy)^{\frac12}S_1\approx 0$ and
$(xy)^{\frac12}S_2\approx 0$. 

The expression $S_3+S_4$ can be split into two parts: the first one is
obtained by replacing in $S_3+S_4$ the functions $A_i(\cdot)$ and
$B_i(\cdot)$ by their constant terms $A_i(0)$ and $B_i(0)$,
respectively, and the second part is a rest term. 

An easy check shows that the first part is equals to 
$$
\multline\frac{1}{\Gamma(z)\Gamma(z')}\Biggl(
\frac{(x/y)^{\frac{z-z'}{2}}
A_1(0)B_2(0)-(x/y)^{\frac{z'-z}{2}}A_1(0)B_2(0)}{x-y}\\
+\frac{(x/y)^{\frac{z'-z}{2}}
B_1(0)A_2(0)-(x/y)^{\frac{z-z'}{2}}B_1(0)A_2(0)}{x-y}\Biggr)\\
=\frac{\sin\pi z\cdot \sin\pi z'}
{\pi\sin\pi(z-z')}\cdot\frac{(x/y)^{\frac{z-z'}{2}}
-(x/y)^{\frac{z'-z}{2}}}{x-y}\,,
\endmultline
$$
where we have used the formula 
$$
\Gamma(z)\Gamma(1-z)=\frac\pi{\sin\pi z}\,.
$$
As for the rest term, it can be written in the form
$$
R(x,y)=\frac{x^\mu y^{-\mu}a(x,y)-x^{-\mu}y^\mu a(y,x)}{x-y}\,,
$$
where $\mu=\frac{z-z'}2$ and $a(x,y)$ is an analytic function near
$(0,0)$ such that $a(0,0)=0$; its exact form is unessential. We shall
prove that $(xy)^{\frac12}R(x,y)\approx 0$. 

Since $a(x,y)$ vanishes at $(0,0)$ we can write it in the form
$$
a(x,y)=x a_1(x,y)+y a_2(x,y)
$$
with certain analytic functions $a_1(x,y)$ and $a_2(x,y)$. With this
notation, we have
$$
(xy)^{\frac12}R(x,y)=\frac{S(x,y)}{x-y}\,,
$$
where
$$
\gather
S(x,y)=x^{\frac32+\mu}y^{\frac12-\mu}a_1(x,y)
+x^{\frac12+\mu}y^{\frac32-\mu}a_2(x,y)\\
-x^{\frac12-\mu}y^{\frac32-\mu}a_1(y,x)
-x^{\frac32-\mu}y^{\frac12+\mu}a_2(y,x).
\endgather
$$

By symmetry, we may assume $x\ge y$. Since the function
$x\mapsto S(x,y)$ vanishes at $x=y$, there exists a point
$\xi\in[y,x]$ such that 
$$
\frac{S(x,y)}{x-y}=\frac{S(x,y)-S(y,y)}{x-y}
=\frac{\partial S(x,y)}{\partial x}\biggm|_{x=\xi}\,.
$$
The same argument as above shows that
$$
x^{\frac12\pm\mu}y^{\frac12\mp\mu}\bigm|_{x=\xi}
=\xi^{\frac12\pm\mu}y^{\frac12\mp\mu}\approx0.
$$
Note also that
$$
\frac yx\biggm|_{x=\xi}=\frac y\xi\,\le\,1.
$$
Looking at $S(x,y)$ we see that $\partial S(x,y)/\partial x$ is the
sum of terms each of which is equal to
$x^{\frac12\pm\mu}y^{\frac12\mp\mu}$ multiplied by an 
analytic function near $(0,0)$ and possibly also multiplied by $y/x$.
Hence, after substitution $x=\xi$, each term is $\approx0$. \qed 
\enddemo

\demo{Proof of Proposition 4.1.4} We shall use the following formula from [E], 6.8(13):
$$
\multline
x^{-1/2}W_{\kappa,0}=\frac{-e^{-x/2}}{\Gamma(1/2-\kappa)}\Bigl(\Phi(1/2-\kappa,1;x)\ln x\\+\sum\limits_{r=0}^{\infty}\frac{(1/2-\kappa)_r}{r!}[\psi(1/2-\kappa+r)-2\psi(1+r)]\frac{x^r}{r!}\Bigr)
\endmultline
$$
where $\psi(x)=\Gamma'(x)/\Gamma(x)$. As in the previous proof, we need only the first terms of the summands. We have
$$
\align
&\varphi_1(x)=C_1(x)\ln x+D_1(x)\\
&\varphi_2(x)=C_2(x)\ln x+D_2(x)
\endalign
$$
where the functions $C_i(x),D_i(x),\  i=1,2,$ are analytic in a neighbourhood of zero and
$$
C_1(0)=-\frac{1}{\Gamma(-z)},\quad D_1(0)=-\frac{-\psi(-z)+2\psi(1)}{\Gamma(-z)};
$$
$$
C_2(0)=-\frac{-1}{\Gamma(1-z)},\quad D_2(0)=-\frac{-\psi(1-z)+2\psi(1)}{\Gamma(1-z)}.
$$
Then by definition of the Whittaker kernel
$$
\gather
K(x,y)=\frac{1}{\Gamma^2(z)}\Biggl(\ln x\ln y\frac{C_1(x)C_2(y)-C_1(y)C_2(x)}{x-y}\\+
\frac{D_1(x)D_2(y)-D_1(y)D_2(x)}{x-y}\\
+\frac{\ln x C_1(x)D_2(y)-\ln y C_1(y)D_2(x)}{x-y}\\+\frac{\ln y D_1(x)C_2(y)-\ln x D_1(y)C_2(x)}{x-y}\Biggr)
\endgather
$$
Let us denote these four summands by $\tilde S_1,\tilde S_2,\tilde S_3,\tilde S_4$, and use the same notation as in the proof of Proposition 4.1.3.  Clearly,
 $$
(xy)^{1/2}|\ln x\ln y|\leq \max\{x \ln^2x,y\ln^2 y\}\approx 0.
$$ 
The functions 
$$
\frac{C_1(x)C_2(y)-C_1(y)C_2(x)}{x-y},\qquad \frac{D_1(x)D_2(y)-D_1(y)D_2(x)}{x-y}
$$
are analytic and, thus, bounded near the origin. Hence,
$$
(xy)^{\frac 12}\tilde S_1\approx(xy)^{\frac 12}\tilde S_2\approx 0.
$$
As in the proof of Proposition 4.1.3, let us split $\tilde S_3+\tilde S_4$ into two parts. In the first one we substitute the functions $C(\cdot)$ and $D(\cdot)$ by their constant terms $C(0)$ and $D(0)$, and the second part is the rest. 
Using the well-known identity
$$
\psi(w)=\psi(1+w)-\frac{1}{w}
$$
 one easily checks that
$$
\multline
\frac{1}{\Gamma^2(z)}\Biggl(\frac{\ln x C_1(0)D_2(0)-\ln y C_1(0)D_2(0)}{x-y}+\frac{\ln y D_1(0)C_2(0)-\ln x D_1(0)C_2(0)}{x-y}\Biggr)\\=\frac{\sin^2\pi z}{\pi^2}\cdot\frac{\ln x-\ln y}{x-y}.
\endmultline
$$
Let us denote the rest term by $\tilde R(x,y)$. We have to show that 
$$
(xy)^\frac{1}{2}\tilde R(x,y)\approx 0.
$$
 We can the rest term in the form
$$
\tilde R(x,y)=\frac{\ln x\cdot b(x,y)-\ln y\cdot b(y,x)}{x-y},
$$
where the function $b(x,y)$ is analytic in a neighbourhood of the origin and $b(0,0)=0$. Then we can represent $b(x,y)$ in the form
$$
b(x,y)=xb_1(x,y)+yb_2(x,y)
$$
where $b_1(x,y)$ and $b_2(x,y)$ are also analytic near the origin. If we introduce the function
$$
\tilde S(x,y)=x^{\frac 32}y^{\frac 12}\ln x b_1(x,y)+x^{\frac 12}y^{\frac 32}\ln x b_2(x,y)-x^{\frac 32}y^{\frac 12}\ln y b_2(y,x) - x^{\frac 12}y^{\frac 32} \ln y b_1(y,x)
$$
then we get
$$
(xy)^\frac{1}{2} \tilde R(x,y)=\frac{\tilde S(x,y)}{x-y}.
$$
 Again, as in the proof of Proposition 4.1.3,  we may assume, by symmetry, that $x\geq y$ and apply the mean value theorem. Thus, for some point $\xi\in [y,x]$
$$
(xy)^\frac{1}{2} \tilde R(x,y)=\frac{\tilde S(x,y)}{x-y}=\frac{\partial \tilde S(x,y)}{\partial x}\biggm|_{x=\xi}\,.
$$
Note that $y/\xi\leq 1$ and
$$
\sqrt{\xi y}\approx 0,\quad \sqrt{\xi y}\ln\xi\approx 0,\quad\frac{y}{\xi}\sqrt{\xi y}\ln\xi\approx 0,
$$
$$
\sqrt{\xi y}\ln y\approx 0,\quad\frac{y}{\xi}\sqrt{\xi y}\ln y\approx 0.
$$
If we explicitly compute the derivative ${\partial \tilde S(x,y)}/{\partial x}$ at $x=\xi$, then we get a sum of expressions estimated above multiplied by analytic in a neighbourhood of the origin functions. Clearly, this sum $\approx 0$.\qed
\enddemo
\enddemo
\subhead 4.2*. Analytic continuation of $F_B^{[m]}(a,b;c|y)$\endsubhead
It turns out that to compute the asymptotics of correlation functions in the non-lifted case is much more difficult than to do that for the lifted processes. The reason is simple. In the lifted case we can use well-known asymptotic formulas for Whittaker functions, while for the original processes we have to deal with asymptotics of multidimensional integrals -- integral representation of Lauricella functions. The asymptotics in both lifted and non-lifted cases happens to be the same. One of possible explanations is that the asymptotic behaviour really depends only on the ratios of the variables (Remark 4.1.2), and lifting does not change them. However, the direct proof of this coincidence is unknown. 

  In this section we present the Mellin-Barnes type integral representation of Lauricella function of type B  and derive certain formulas for analytic continuation of  $F_B^{[m]}(a,b;c|y)$. These expansions  will be used for computing the asymptotics.

The results of this section represent a generalization of the well known formulas for the analytic continuation of the Gauss hypergeometric function and Appell hypergeometric function $F_3$, see [E], [Mar], [Ex1], [Ex2].

We start with Mellin-Barnes type integrals.
\proclaim{Proposition 4.2.1} If $a_i,b_i\ne 0,-1,-2,\ldots$ for all $i=1,\ldots,m$ then
$$
\multline
F_B^{[m]}(a,b;c|y) =  \frac{\Gamma(c)}{\prod\limits_{i=1}^m\Gamma(a_i)\Gamma(b_i)}\\ \times
\frac{1}{{(2\pi i)}^m}\int\limits_{-i\infty}^{+i\infty}\cdots \int\limits_{-i\infty}^{+i\infty}
\frac{\prod\limits_{i=1}^m \Gamma(a_i+s_i)\Gamma(b_i+s_i)\Gamma(-s_i){(-y_i)}^{s_i}ds_i}{\Gamma\left(c+\sum\limits_{i=1}^ms_i\right)}
\endmultline
$$
where $\arg{(-y_i)}<\pi$ and the $i$th path of integration separates the points $s_i=0,1,2,\ldots$ from the points $s_i=-a_i-n,\ s_i=-b_i-n\ (n=0,1,\ldots).$
\endproclaim
\demo{Sketch of the proof}
This formula can be found in [Ex1, (2.5.6)] and [Ex2, (5.2.3.7)]. One can  obtain the proof by  computing the residues of the integrand at the points  $s_i=0,1,2,\ldots$ for all $i=1,\ldots,m$. Then one gets exactly the series (2.7). The correctness of this operation can be checked by the general techniques described in [Mar]. \qed
\enddemo
Note that if we want to analyze the behaviour of $\rz$ at the origin, then, see Theorem 2.4.1 and (2.9), we need to know the behaviour of $F_B^{[m]}(a,b;c|y)$ for large negative values of variables.
The Mellin-Barnes  integral representation allows us to continue the Lauricella function to this domain.
\proclaim{Proposition 4.2.2}
Let $\Re y_1,\ldots,\Re y_m$ be negative and sufficiently large; and $a_i-b_i\notin \Bbb Z$ for all $i=1,\ldots,m$. Then
$$
\align
F_B^{[m]}(a,b;c|y)=&\frac{\Gamma(c)}{\prod\limits_{i=1}^m\Gamma(a_i)\Gamma(b_i)}
\sum\limits_{\Sb \Cal I =\{i_1<\ldots<i_p\}\subset\{1,\ldots,m\}\\
\overline{\Cal I}=\{j_1<\ldots<j_{m-p}\}\endSb} \prod\limits_{k=1}^p{(-y_{i_k})}^{-a_{i_k}}
\prod\limits_{l=1}^{m-p}{(-y_{j_l})}^{-b_{j_l}}\\
\times\sum\limits_{\Sb \alpha_1,\ldots,\alpha_p\geq 0\\ \beta_1,\ldots,\beta_{m-p}\geq 0\endSb}
&\prod\limits_{k=1}^p
\frac{\Gamma(b_{i_k}-a_{i_k}-\alpha_k)\Gamma(a_{i_k}+\alpha_k)}
{\alpha_k!{(y_{i_k})}^{\alpha_k}}
\prod\limits_{l=1}^{m-p}\frac
{\Gamma(a_{j_l}-b_{j_l}-\beta_l)\Gamma(b_{j_l}+\beta_l)}{\beta_l!{(y_{j_l})}^{\beta_l}}\\
&\qquad\qquad\qquad\qquad\qquad\times\frac{1}{\Gamma\left(c-\sum\limits_{k=1}^p(a_{i_k}+\alpha_k)-\sum\limits_{l=1}^{m-p}(b_{j_l}+\beta_l)\right)}.
\tag 4.1
\endalign
$$
\endproclaim
\example {Remark 4.2.3} This formula represents $F_B^{[m]}(a,b;c|y)$ as a finite sum of
expressions of the form `product of minus variables in some (complex) powers times a function analytic at the infinity'. For $m=1$ this is the well-known formula for analytic continuation of the Gauss hypergeometric function, see [E], 2.1.4(17):
$$
\multline
\frac{F(a,b;c;w)}{\Gamma(c)}=\frac{\Gamma(b-a)}{\Gamma(b)\Gamma(c-a)}\frac{1}{(-w)^a}F(a,1-c+a;1-b+a,w^{-1})\\+\frac{\Gamma(a-b)}{\Gamma(a)\Gamma(c-b)}\frac{1}{(-w)^b}F(b,1-c+b;1-a+b,w^{-1}).
\endmultline
$$
For $m=2$ the formula consists of four summands and is a known example when one hypergeometric series ($F_3$ in this case) is continued by another (it will be $F_2$ here), see [E], 5.11(10).
\endexample
\demo {Sketch of the proof} The formula can be easily obtained by calculating the residues of the integrand from the Mellin-Barnes type integral in Proposition 4.2.1 at the points  $s_i=-a_i-k,\ s_i=-b_i-k\ (k=0,1,\ldots).$ The correctness of this operation, as in Proposition 4.2.1, is checked by the machinery developed in [Mar]. The set $ \Cal I\subset\{1,\ldots,m\} $ reflects the fact that for every $i\in \{1,\ldots,m\}$ we have two sequences of poles, if we use the first sequence (i.e.,
$s_i=-a_i-k$ for some $k=0,1,\ldots$) then we put $i$ in $\Cal I$; if we use the second sequence
(i.e., $s_i=-b_i-k$ for some $k=0,1,\ldots$) then $i\in \overline{\Cal I}= \{1,\ldots,m\}\setminus \Cal I$.\qed
\enddemo
We shall apply the last proposition to the hypergeometric functions used in Chapter 2. The
restriction $a_i-b_i\notin \Bbb Z$ implies for us that $z-z'\notin \Bbb Z$ (see (2.9')) which is
equivalent to $z\ne z'$. In case $z=z'$ we have $a_i=b_i$ for all $i=1,\ldots,m$, and all the
poles $s_i=-a_i-n,\  (n=0,1,\ldots)$ of the integrand in Proposition 4.2.1 are of the second
order. This `logarithmic' case can be worked out by the same techniques as the ordinary one
(i.e., $z\neq z'$). Below we present the `logarithmic' version of Proposition 4.2.2. Its special case for $m=1$ is the known formula  [E], 2.1.4(18). 
\proclaim{Proposition 4.2.4} 
 $$
\multline
F_B^{[m]}(a,a;c|y)=\frac{\Gamma(c)}{\prod\limits_{i=1}^m\Gamma^2(a_i)}
\su_{k_1,\ldots,k_m}\pro_{i=1}^m (-y_i)^{-a_i}\frac{\Gamma(a_i+k)}{(k!)^2}\frac{(\operatorname{ln}(-y_i)+h_{k_i}(a))}{(-y_i)^{k_i}}\\
\times \frac{1}{\Gamma\left(c-\su_{i=1}^m(a_i+k_i)\right)}
\endmultline
$$
where
$$
h_{k_i}(a)=2\psi(k_i+1)-\psi(a_i+k_i)-\psi\left(c-\su_{i=1}^m(a_i+k_i)\right),\quad \psi(x)=\frac{\Gamma'(x)}{\Gamma(x)}.
$$
\endproclaim

\subhead 4.3*. Asymptotics of the non-lifted correlation functions \endsubhead 
In this section we shall prove the following result, cf. Theorem 4.1.1.
\proclaim{Theorem 4.3.1} Let $x_1,\ldots,x_n>0$. Then
$$
\rho_n^{(zz')}(x_1,\ldots,x_n)=\frac{\det k(x_i/x_j)
+\tilde r(x_1,\ldots,x_n)}{x_1\cdots x_n}
$$
where $k(x)$ is defined in Theorem 4.1.1 and
$$
\tilde r(x_1,\ldots,x_n)=\cases
O\left(\max\{x_i\}\right),\ &\text{if}\ z'=\bar z\ne z\\
O\left(\max\{x_i\ln^2 x_i\}\right),\ &\text{if}\ z'=z\in \Bbb C\setminus \Bbb Z\\
O\left({(\max\{x_i\})}^{1-|z-z'|}\right),\ &\text{if}\ m<z,z'<m+1,\ m\in\Bbb Z
\endcases
$$
as $x_1,\ldots,x_n\to 0$.
\endproclaim 
\example{Remark 4.3.2} The asymptotics coincides with that in the lifted case.
\endexample
\demo{Proof of Theorem 4.3.1} We shall use the notation from Section 2.4. As we have seen (Theorem 2.4.1),
$$
\align 
\rz(x_1,\ldots,x_n)=\Gamma(t)\prod\limits_{i=1}^n&\phi_{z-1}(x_i)\phi_{z'-1}(x_i)\cdot\phi_{c-1}(1-|x|)\\
\times &\sum\limits_{\sigma\in S_n}\operatorname{sgn} \sigma \cdot f_n^{(zz')}(y_1,\ldots,y_n;y_{\sigma(1)},\ldots,y_{\sigma(n)}).
\endalign
$$
where $c=t-n(z+z'-1)$ and for all $i=1,\ldots,n$
$$
y_i=-\frac{1-|x|}{x_i},\quad |x|=x_1+\ldots+x_n,
$$
 the function $f$ is defined in (2.9).

Thus, it suffices to show that
$$
\multline
\Gamma(t)\prod\limits_{i=1}^n\phi_{z-1}(x_i')\phi_{z'-1}(x_i'')\cdot\phi_{c-1}\left(1-\frac{|x'|+|x''|}{2}\right)f_n^{(zz')}(y_1',\ldots,y_n';y_{1}'',\ldots,y_{n}'')\\ = \frac{\pro_{i=1}^n
 k({x_i'}/{x_i''})+\hat r(x_1',\ldots,x_n';x_1'',\ldots,x_n'')}{(x_1'\cdots x_n')^{\frac 12 +\frac{z'-z}{2}}(x_1''\cdots x_n'')^{\frac 12 +\frac{z-z'}{2}}}
\endmultline
$$
where
$$
y_i'=-\frac{1-\frac{ |x'|+|x''|}{2}}{x_i'},\quad y_i''=-\frac{1-\frac{ |x'|+|x''|}{2}}{x_i''}
$$
and
$$
\hat r(x',x'')=\cases
O\left(\max\{x_i',\ x_i''\}\right),\ &\text{if}\ z'=\bar z\ne z\\
O\left(\max\{x_i'\ln^2 x_i',\  x_i''\ln^2 x_i'' \}\right),\ &\text{if}\ z'=z\in \Bbb C\setminus \Bbb Z\\
O\left({(\max\{x_i',\ x_i''\})}^{1-|z-z'|}\right),\ &\text{if}\   m<z,z'<m+1,\ m\in\Bbb Z
\endcases
$$
as all $x_i',x_i''$ tend to zero.

We shall give the proof for $n=1$, for $n>1$ the proof is obtained by applying the same arguments coordinate-wise. 

As for the lifted processes, the cases $z\neq z'$ and $z=z'$ are different. First, let $z\neq z'$.
Then by (2.9)
$$
f(y',y'')=\frac{1}{y'-y''}\left(y'F_B^{[2]}(a^0,b^0;c|y)-y''F_B^{[2]}(a^1,b^1;c|y)\right)
\tag 4.2
$$
where $c=(1-z)(1-z')$;
$$
a^0=(1-z';-z),\quad a^1=(-z';1-z);$$$$ b^0=(1-z;-z'),\quad b^1=(-z;1-z');
$$
$$
y=(y';y'').
$$
Now we apply Proposition 4.2.2 (or [E], 5.11(10)) to two hypergeometric functions in (4.2). We get
$$
\gather
F_B^{[2]}(a^0,b^0;c|y)=\frac{\Gamma(c)}{\Gamma(t)}\Biggl((-y')^{z'-1}(-y'')^{z}\frac{\Gamma(z'-z)\Gamma(z-z')}{\Gamma(1-z)\Gamma(-z')}G_1(1/y',1/y'')\\
+(-y')^{z'-1}(-y'')^{z'}\frac{\Gamma(z'-z)\Gamma(z'-z)}{\Gamma(-z)\Gamma(1-z)}G_2(1/y',1/y'')
\\
+(-y')^{z-1}(-y'')^{z}\frac{\Gamma(z-z')\Gamma(z-z')}{\Gamma(-z')\Gamma(1-z')}G_3(1/y',1/y'')\\
+(-y')^{z-1}(-y'')^{z'}\frac{\Gamma(z'-z)\Gamma(z-z')}{\Gamma(1-z')\Gamma(-z)}G_4(1/y',1/y'')\Biggr)
\endgather
$$
where $G_i$'s are analytic at the origin and $G_i(0,0)=1$ for all $i=1,2,3,4$.
Similarly,
$$
\gather
F_B^{[2]}(a^0,b^0;c|y)=\frac{\Gamma(c)}{\Gamma(t)}\Biggl((-y')^{z'}(-y'')^{z-1}\frac{\Gamma(z'-z)\Gamma(z-z')}{\Gamma(1-z')\Gamma(-z)}\hat G_1(1/y',1/y'')\\
+(-y')^{z'}(-y'')^{z'-1}\frac{\Gamma(z-z')\Gamma(z-z')}{\Gamma(-z)\Gamma(1-z)}\hat G_2(1/y',1/y'')
\\
+(-y')^{z}(-y'')^{z-1}\frac{\Gamma(z-z')\Gamma(z-z')}{\Gamma(-z')\Gamma(1-z')}\hat G_3(1/y',1/y'')\\
+(-y')^{z}(-y'')^{z'-1}\frac{\Gamma(z'-z)\Gamma(z-z')}{\Gamma(1-z)\Gamma(-z')}\hat G_4(1/y',1/y'')\Biggr)
\endgather
$$
where, again, $\hat G_i$'s are analytic at the origin and $\hat G_i(0,0)=1$ for all $i=1,2,3,4$.
After we plug  these expressions into (4.2) and rewrite the result in terms of solely $x_i'$'s and $x_i''$'s, we shall get the expression  of the same form (modulo multiplication by $(x/y)^{(z-z')/2}$), as we had for $K(x,y)$ in the proof of Proposition 4.1.3. Then we word for word follow this proof, and, thus, prove the assertion for $z\neq z'$.

If $z=z'$ then (4.2) still holds, but we need to use different expressions for the hypergeometric function, namely, we use Proposition 4.2.4. We get
$$
\gather
F_B^{[2]}(a^0,b^0;c|y)=\frac{\Gamma(c)(-y')^{z-1}(-y'')^z}{\Gamma(t)\Gamma(1-z)\Gamma(-z)}\Bigl({\ln(-y')\ln(-y'')}L_1(1/y',1/y'')\\
+{\ln(-y')}(2\psi(1)-\psi(-z)-\psi(t))L_2(1/y',1/y'')\\+
\ln(-y'')(2\psi(1)-\psi(1-z)-\psi(t))L_3(1/y',1/y'')\\+
(2\psi(1)-\psi(-z)-\psi(t))(2\psi(1)-\psi(1-z)-\psi(t))L_4(1/y',1/y'')
\Bigr)
\endgather
$$
where $L_i$'s are analytic at the origin, $L_i(0,0)=1$ for all i;
and
$$
\gather
F_B^{[2]}(a^1,b^1;c|y)=\frac{\Gamma(c)(-y')^{z}(-y'')^{z-1}}{\Gamma(t)\Gamma(1-z)\Gamma(-z)}\Bigl({\ln(-y')\ln(-y'')}\hat L_1(1/y',1/y'')\\
+{\ln(-y')}(2\psi(1)-\psi(1-z)-\psi(t))\hat L_2(1/y',1/y'')\\+
\ln(-y'')(2\psi(1)-\psi(-z)-\psi(t))\hat L_3(1/y',1/y'')\\+
(2\psi(1)-\psi(-z)-\psi(t))(2\psi(1)-\psi(1-z)-\psi(t))\hat L_4(1/y',1/y'')
\Bigr)
\endgather
$$
where, as usual $\hat L_i$'s are analytic at the origin and their constant terms are equal to one. 

Again, after we plug  this expressions in (4.2) we shall get an expression which will coincide with the expression for $K(x,y)$ obtained in the proof of Proposition 4.1.4. We deal with it exactly as we did in that proof, and, finally, prove the theorem. \qed
\enddemo

\enddocument